\documentclass{article}

\usepackage[utf8]{inputenc}
\usepackage{amsfonts}
\usepackage{amsthm}
\usepackage{amsmath}
\usepackage{xcolor}
\usepackage{comment}
\usepackage{enumerate}
\usepackage{authblk}

\usepackage[style=alphabetic]{biblatex}
\newtheorem{Thm}{Theorem}[section]
\newtheorem{Lem}[Thm]{Lemma}
\newtheorem{Def}[Thm]{Definition}
\newtheorem{Prop}[Thm]{Proposition}
\newtheorem{Cor}[Thm]{Corollary}
\newtheorem*{Rem}{Remark}

\newcommand{\mb}{\mathbb}
\newcommand{\mc}{\mathcal}

\numberwithin{equation}{section}

\title{The Yamabe flow on asymptotically flat manifolds}
\date{\vspace{-5ex}}
\author{Eric Chen\thanks{University of California, Santa Barbara, ecchen@math.ucsb.edu; research partially supported by an AMS--Simons Travel grant} \qquad Yi Wang\thanks{Johns Hopkins University, ywang261@jhu.edu;
research partially supported by NSF CAREER Award DMS-1845033}}

\begin{document}
\maketitle

\begin{abstract}
We study the Yamabe flow starting from an asymptotically flat manifold $(M^n,g_0)$.
%of order $\tau>\frac{n-2}{2}$. 
We show that the flow converges to an asymptotically flat, scalar flat metric in a weighted global sense if $Y(M,[g_0])>0$, and show that the flow does not converge otherwise. If the scalar curvature is nonnegative and integrable, then the ADM mass at time infinity drops by the limit of the total scalar curvature along the flow.

%As a consequence of our weighted convergence results we find that along a Yamabe flow starting from an asymptotically flat manifold with nonnegative, integrable scalar curvature the ADM mass weakly decreases along the flow out to time infinity, and then possibly drops at the limit.
\end{abstract}

\section{Introduction}

In this article we study the long-time existence and convergence of the Yamabe flow
\begin{align}
\begin{cases}
\frac{\partial g}{\partial t}=-R g,\label{floweq}
\\
g(0)=g_0,
\end{cases}
\end{align}
starting from an asymptotically flat manifold $(M^n,g_0)$. Here $R$ denotes the scalar curvature of the Riemannian metric $g=g(t)$. This flow preserves the conformal class of $g$ and is the natural analogue of the volume-normalized Yamabe flow on compact manifolds introduced by Hamilton \cite{Hamilton}. He was motivated by the resolution of the Yamabe problem due to Yamabe, Trudinger, Aubin, and Schoen \cite{Yamabe,Trudinger,Aubin,Schoen}, which showed that every conformal class of Riemannian metrics on a compact manifold admits at least one metric of constant scalar curvature, known as a Yamabe metric, which minimizes the Einstein--Hilbert functional. Hamilton proposed the volume-normalized Yamabe flow, which can be viewed as the gradient flow of the Einstein--Hilbert functional within a fixed conformal class, as a natural evolution equation which could potentially evolve a given metric on a compact manifold to one of constant scalar curvature within the same conformal class.

Hamilton was already able to prove %in general 
the long-time existence of the volume-normalized Yamabe flow on compact manifolds when he introduced it. Convergence of the flow to a metric of constant scalar curvature has now been mostly settled by work of Ye, Schwetlick--Struwe, and Brendle \cite{Ye,Struwe,Brendle5,Brendle6}. 
%In \cite{Brendle6} one still requires a technical assumption on the vanishing of the Weyl tensor in dimensions six and higher.

The study on noncompact manifolds is less developed, but there have been a number of long-time existence results for the flow \eqref{floweq}. For complete noncompact manifolds, Ma has shown long-time existence when starting from a metric of non-negative scalar curvature which is conformal to a metric of non-positive scalar curvature, assuming the corresponding conformal factor is bounded from above \cite{MaLi}. Schulz has shown long-time existence existence under related hypotheses --- namely when starting from a metric of positive Yamabe constant which is conformal to a metric of non-positive and bounded scalar curvature, assuming the corresponding conformal factor is bounded both from above and from below \cite{Schulz}. Other works give long-time existence results in the settings of conformally hyperbolic and singular spaces \cite{Schulz2,Bahuaud,Olsen}. Similar to the compact case, convergence results for \eqref{floweq} have been slower to develop --- we are aware only of Ma's work mentioned above \cite{MaLi} in which he is also able to show $C^\infty_{\text{loc}}$ convergence to a scalar flat limit metric, using crucially the assumption of initially non-negative scalar curvature.
%, along with a few more recent results which we mention at the end of this section.

Our work adds to the study of \eqref{floweq} in the non-compact case by providing a long-time existence result in the setting of asymptotically flat manifolds without requiring any additional curvature assumptions, and by showing that the flow in this setting converges in a strong, global weighted sense (in $C^\infty_{-\tau'}$) to a scalar flat, asymptotically flat metric whenever one might hope for this --- namely whenever there exists a scalar flat, asymptotically flat metric lying in the conformal class of the initial metric. Asymptotically flat metrics are interesting to study under Yamabe flow because as shown by  Cheng--Zhu, who were motivated by similar results in the Ricci flow setting \cite{DaiMa}, asymptotic flatness is preserved under the Yamabe flow and moreover the ADM mass is monotonically decreasing \cite{CZ}.  As a consequence of our convergence results, the drop from the ADM mass along the flow to the mass of the limiting scalar flat metric is accounted for by the total scalar curvature pushed out to spatial infinity by the flow.

\begin{comment}

A week before we completed this manuscript we learned of recent new work by Ma in some directions similar to ours. He proves a version of the long-time existence of the flow starting from any asymptotically flat manifold as in our Theorem \ref{ltegen} below, and further claims the uniqueness of this flow \cite{LiMa21}. He additionally has a result on the $C^\infty_{loc}$ convergence of the flow to a scalar flat metric when starting with $R_{g_0}\geq 0$ which appears to be similar to the aforementioned convergence result in \cite{MaLi}, and is a consequence of a more general $C^\infty_{loc}$ convergence result for metrics in the conformal class of a scalar nonnegative metric satisfying a technical condition \cite{LiMa20}. In comparison we establish a uniform decay rate estimate of the scalar curvature with which we obtain global weighted and unweighted convergence in the more general case when $Y(M^n,[g_0])>0$. We also exclude the possibility of convergence in the remaining case $Y(M^n,[g_0])\leq 0$.

\end{comment}

While we were completing this manuscript we learned of recent new work by Ma \cite{LiMa21} in some directions similar to ours. His result considers the long-time existence of the flow, which overlaps with some of the content in our Lemma \ref{utimebd} and its corollary Theorem \ref{ltegen}. But the focus of our paper is to prove the strong global $C^\infty_{-\tau'}$ convergence when $Y(M^n,[g_0])>0$. Towards this goal, the central issue is the sufficiently fast decay rate estimate of $\|R\|_{L^\infty}$, which we achieve in several steps. One of the first observations is the monotonicity of the integral $\|R\|_{L^{p}}$ for $p= n/2$. This allows us to prove that $\|R\|_{L^{p}}$ tends to zero as $t\rightarrow \infty$ for $p$ in a neighborhood of $n/2$. Another key step is to obtain the convergence of $\|R\|_{L^{\frac{np}{n-2}}}\rightarrow 0$ at a particular decay rate 
%$o(t^{-\frac{n}{(n-2)q}})$ for suitable $L^q$ norms of $R$. 
using the Moser iteration. With these, we finally conclude $\|R\|_{L^{\infty}}\leq O(t^{-1-\delta})$ --- here $\delta>0$ is essential as the integrability of $\int_0^\infty\|R\|_{L^{\infty}}\ dt$ leads to the unweighted and weighted convergence of the flow in Section 5.
To better understand the complete picture, we also study the flow when $Y(M^n,[g_0])\leq 0$.
In this remaining case, we exclude the possibility of convergence.

\begin{comment}

We take the long-time existence as a starting point. Then We observed the monotonicity of the integral bound of the scalar curvature implies the scalar curvature tends to $0$. We then improve such an estimate to the uniform decay estimate of  the integral bound of the scalar curvature, and finally get the uniform decay estimate of the
$\|R\|_{L^\infty}$ 
could be improved to the $L^\infty$ decay estimate
to obtain global weighted and unweighted convergence in the more general case when $Y(M^n,[g_0])>0$.

starting from any asymptotically flat manifold as in our Theorem \ref{ltegen} below, and further claims the uniqueness of this flow \cite{LiMa21}. He additionally has a result on the $C^\infty_{loc}$ convergence of the flow to a scalar flat metric when starting with $R_{g_0}\geq 0$ which appears to be similar to the aforementioned convergence result in \cite{MaLi}, and is a consequence of a more general $C^\infty_{loc}$ convergence result for metrics in the conformal class of a scalar nonnegative metric satisfying a technical condition \cite{LiMa20}. In comparison we establish a uniform decay rate estimate of the scalar curvature with which we obtain global weighted and unweighted convergence in the more general case when $Y(M^n,[g_0])>0$. We also exclude the possibility of convergence in the remaining case $Y(M^n,[g_0])\leq 0$.

\end{comment}

\subsection{Main results}

Below we describe our main results. The definition of a $C^{k+\alpha}_{-\tau}$ asymptotically flat (AF) manifold below is the same as in \cite{CZ,MaLi} and this along with related notions of weighted H\"{o}lder and Sobolev spaces are stated precisely in Section \ref{analyticprelim}.
%By asymptotic flatness we mean the following:
%\begin{Def}\label{defaf}
%We call a smooth Riemannian manifold $(M^n,g)$ asymptotically flat of order $\tau>0$ if for some compact set $K\subset M^n$, there exists an $R_0>0$ and a diffeomorphism $\Phi:M^n\backslash K\rightarrow \mb{R}^n\backslash B_{R_0}(0)$ such that for some $\tau>0$,
%\begin{align}\notag
%g_{ij}(x)=\delta_{ij}+O(|x|^{-\tau})\quad\text{and}\quad\partial^\alpha g_{ij}(x)=O(|x|^{-\tau+|\alpha|}),
%\end{align}
%for partial derivatives $\partial^\alpha$ of any order, as $|x|\rightarrow\infty$ in $\mb{R}^n$. $\Phi$ is called the asymptotically flat coordinate system, and we write $M_\infty=M\backslash K$.
%\end{Def}

%Note that this implies $\int |R|^p\ dV_t<\infty$ for all $p>\frac{n}{2+\tau}$.
Our first result is the long-time existence of the Yamabe flow on all $AF$ manifolds. On the order of asymptotic flatness we always assume that $\tau>0$. 
\begin{Thm}\label{ltegen}
Given any $C^{k+\alpha}_{-\tau}$ AF manifold $(M^n,g_0)$, $k\geq 2$, there exists a Yamabe flow starting from it defined for all positive times with $(M^n,g(t))$ remaining $C^{k+\alpha}_{-\tau'}$ asymptotically flat for all $\tau'\leq\min\{\tau,n-2\}$.
\end{Thm}

We will not be concerned in this work with the uniqueness of Yamabe flow in the AF setting; above, Theorem \ref{ltegen} refers to any fine solution of the Yamabe flow whose short-time existence on AF manifolds was shown in \cite[Corollary 2.5]{CZ} --- see Definition \ref{finedef} for details. These are always the solutions which we study in this work, but below we will often write \emph{the Yamabe flow} with this meaning implicit for ease of presentation.

The special case $R_{g_0}\geq 0$ of Theorem \ref{ltegen} was previously known from work of Ma \cite[Theorem 1]{MaLi}, while work of Schulz implies the case  $Y(M,[g_0])>0$ \cite[Theorem 1]{Schulz}. As mentioned earlier, while we were preparing this manuscript we also learned of recent work of Ma \cite{LiMa21} establishing the $C^{2+\alpha}_{-\tau}$ version of the long-time existence result of Theorem \ref{ltegen}.
%in which he also claims uniqueness of the flow \cite{LiMa21}.

Here, $Y(M,[g_0])$ denotes the following conformally invariant quantity which we call the Yamabe constant, motivated by the definition of the Yamabe constant in the compact case:
\begin{align}
Y(M,[g_0]):=\inf_{\substack{v\in C^{\infty}_0(M),\\{v\neq 0}}} \frac{4\frac{n-1}{n-2}\int_M|\nabla v|^2+R_{g_0} v^2\ dV_{g_0}}{\left(\int |v|^{\frac{2n}{n-2}}\ dV_{g_0}\right)^{\frac{n-2}{n}}},
\end{align}
which plays an important role in the prescribed scalar curvature problem on conformal classes of AF metrics \cite{CB,Maxwell,DM}. 
%We write the condition $Y(M,[g_0])>0$ above because 
For instance, by \cite[Proposition 3]{Maxwell}, which gives the correct version of a result first claimed in \cite[Theorem 2.1]{CB}, the conformal class $(M,[g_0])$ admits a scalar flat, AF metric if and only if it is Yamabe-positive --- we state a version of this result in Proposition \ref{CBY} below.

Having obtained long-time existence of the Yamabe flow on AF manifolds, we next study its convergence properties. On compact manifolds, proving the convergence of the Yamabe flow required much more work than proving long-time existence --- see
%, and in fact is not yet completely settled. 
%We refer interested readers to see the aforementioned work of Ye, Schwetlick--Struwe, and Brendle
\cite{Ye,Struwe,Brendle5,Brendle6}. In our setting, if the conformal class $(M,[g_0])$ admits a scalar flat, asymptotically flat metric (which must be unique), this gives additional information with which we are able to obtain strong quantitative decay of the scalar curvature. By our earlier discussion, we can equivalently formulate such a condition in terms of the positivity of $Y(M,[g_0])$, and we have two distinct possibilities for the behavior of the flow as $t\rightarrow\infty$ as below.

\begin{Thm}\label{uconverge}
Let $(M^n,g_0)$ be a $C^{k+\alpha}_{-\tau}$ AF manifold with $k\geq 3$.
\begin{enumerate}[(1)]
\item\label{uconverge1} If $Y(M^n,[g_0])>0$, then the Yamabe flow $(M^n,g(t))$ starting from $(M^n,g_0)$ converges uniformly in $C^{k+\alpha}_0$ to the unique $C^{k+\alpha}_{-\tau}$ AF metric $g_\infty\in[g_0]$ as $t\rightarrow\infty$.
\item\label{uconverge2} If $Y(M^n,[g_0])\leq 0$, then then the Yamabe flow $(M^n,g(t))$ starting from $(M^n,g_0)$ does not converge. In particular, $g(t)=u(t)^{\frac{4}{n-2}}g_0$ will fail to remain uniformly equivalent to $g_0$ as $t\rightarrow\infty$, and both  $\|u(t)\|_{L^\infty}$ and the $L^2$ Euclidean-type Sobolev constant of $g(t)$ will tend to positive infinity.
\end{enumerate}
\end{Thm}

The second part of Theorem \ref{uconverge} raises the question of what more can be determined regarding the non-convergence along the Yamabe flow on a large class of AF manifolds. 
In the compact case, Schwetlick--Struwe showed that a failure of convergence must imply a particular kind of infinite time bubbling behavior for the volume-normalized Yamabe flow \cite{Struwe}.
%, and this has now been ruled out in nearly all cases now by Brendle \cite{Brendle5,Brendle6}. 
In the non-compact, conformally flat setting Choi--Daskalopoulos have produced examples of Yamabe flows with infinite-time Type II singularities, which satisfy $\sup_{M\times[0,\infty)}|\text{Rm}(x,t)|=\infty$ \cite{ChoiDaskalopoulos}. Finite-time singularities in the non-compact, conformally flat setting have also been studied in \cite{CDK,DKS}.

The first part of Theorem \ref{uconverge} above gives a uniform global convergence which is strong enough to allow us to identify the limiting metric $g_\infty$. If we impose some mild restrictions on the decay order $\tau$ (which are natural for instance if we wish to consider the ADM mass under the flow), then we can further improve to weighted convergence in Theorems \ref{Ypos} and \ref{Ypos3} below.

\begin{Thm}\label{Ypos}
Let $(M^n,g_0)$ be a $C^{k+\alpha}_{-\tau}$ AF manifold with $Y(M,[g_0])>0$, $k\geq 3$, and  $\tau>1$. Then there exists a Yamabe flow $(M^n,g(t))$ starting from $(M^n,g_0)$ defined for all positive times and a metric $g_\infty$ on $M^n$ which is $C^{k+\alpha}_{-\tau'}$ AF for all $\tau'<\min\{\tau,n-2\}$ so that for any such $\tau'$ we have
\begin{align}
\|g(t) - g_\infty \|_{ C^{k+\alpha}_{-\tau'} }=O(t^{- \delta_0  }),  \quad  \mbox{as} \ t\rightarrow \infty,
\end{align}
for some $\delta_0>0$. In particular, this Yamabe flow converges in $C^{k+\alpha}_{-\tau'}$ to the asymptotically flat, scalar flat metric $g_\infty$.
\end{Thm}

%\begin{Rem}
%If we instead start with a
%\end{Rem}

As noted earlier, the work of
%https://www.overleaf.com/project/5f861b9a1418a800016ab4a0
Cheng--Zhu shows that under appropriate initial conditions, the ADM mass of an asymptotically flat manifold with non-negative and integrable scalar curvature is nonincreasing under Yamabe flow \cite[Theorem 1.5]{CZ}, which suggests further study of the mass along the flow. Theorem \ref{Ypos} will suffice for this purpose when $n\geq 4$, since mass is well-defined for $\tau>\frac{n-2}{2}\geq 1$; however, when $n=3$ the condition $\tau>1$ forces the mass to vanish and is too restrictive.

However, by adding the natural conditions $R_{g_0}\geq 0$ and $R_{g_0}\in L^1$ (if we are concerned with the mass) we can still obtain a weighted convergence result that will allow us to study the mass in dimension $n=3$.

\begin{Thm}\label{Ypos3}
Let $(M^3,g_0)$ be a $C^{k+\alpha}_{-\tau}$ AF manifold with $Y(M,[g_0])>0$, $k\geq 3$, $\tau>\frac{1}{2}$, $R_{g_0}\geq 0$, and $R_{g_0}\in L^1$. Then there exists a Yamabe flow $(M^n,g(t))$ starting from $(M^n,g_0)$ defined for all positive times and a metric $g_\infty$ on $M^n$ which is $C^{k+\alpha}_{-\tau'}$ AF for all $\tau'<\min\{\tau,1\}$ so that for any such $\tau'$ we have
\begin{align}
\|g(t) - g_\infty \|_{ C^{k+\alpha}_{-\tau'} }=O(t^{- \delta_0  }),  \quad  \mbox{as} \ t\rightarrow \infty,
\end{align}
for some $\delta_0>0$. In particular, this Yamabe flow converges in $C^{k+\alpha}_{-\tau'}$ to the asymptotically flat, scalar flat metric $g_\infty$.
\end{Thm}

It is well-known that $m(g_0)\geq m(g_\infty)$, since a scalar flat AF metric minimizes the mass among scalar nonnegative metrics within its conformal class \cite{SY}, and one can moreover compute this difference as a multiple of the integral of $R_{g_0}$ against $u_\infty$ with respect to the volume form of $g_0$,
%\begin{align}
%    m(g_0)-m(g_\infty)=\frac{n-2}{4(n-1)\omega_{n-1}}\int R_{g_0}u_\infty\ dV_{g_0},\notag
%\end{align} 
where $g_\infty=u_\infty^{\frac{4}{n-2}}g_0$. By our weighted convergence results in Theorems \ref{Ypos} and \ref{Ypos3} combined with the monotonicity of mass along Yamabe flow \cite{CZ} and the lower semicontinuity of the ADM mass under $C^{2}_{-\tau}$ convergence when $\tau>\frac{n-2}{2}$ \cites{YuLi}[Theorem 14]{McFeron}[Theorem 13]{Jauregui}, this difference is also the mass drop at time infinity of the Yamabe flow starting from $(M^n,g_0)$, and is therefore controlled by the $L^1$ norm of the scalar curvature as it escapes to infinity.

% is to compare the mass along the flow to the mass of the limiting metric under Yamabe flow, which we know exists by Theorem \ref{uconverge}. Motivated by this we prove two weighted convergence results which improve upon the convergence in the first part of Theorem \ref{uconverge} above.

%Although \cite[Theorem 1]{MaLi} shows that the flow converges in $C^\infty_{\text{loc}}$ to a scalar flat metric, one does not know if the limit is AF. Theorem \ref{Ypos} shows that the limit is indeed AF and moreover that we have $C^{k+\alpha}_{-\tau'}$ convergence. 

%when $n\geq 4$.

%Note that the ADM mass is defined for $\tau>\frac{1}{2}=\frac{n-2}{2}$, so our restrictions on $\tau$ above allow us to study the ADM mass along the flow whenever it exists, for all dimensions $n\geq 3$.
%It is possible that $\tau>1$ above is a technical restriction, but at least when $n\geq 4$, so that $\frac{n-2}{2}\geq 1$, this still allows us to mass along the flow.

\begin{Cor}\label{masscor}
For $n\geq 3$, let $(M^n,g_0)$ be a $C^{k+\alpha}_{-\tau}$ AF manifold with non-negative scalar curvature and $k\geq 3$, $\tau>\frac{n-2}{2}$, along with $R_{g_0}\in L^1(M^n,g_0)$. Then along the Yamabe flow $(M^n,g(t))$ starting from $(M^n,g_0)$,
\begin{align}
    \lim_{t\rightarrow\infty}\left(m(g(t))-\frac{1}{2(n-1)\omega_{n-1}}\int R_{g(t)}\ dV_t\right)=m(g_\infty).
\end{align}
In particular if $n=3,4,$ or $5$ or if $\tau>n-3$ so that $m(g(t))$ is constant along the flow, then
\begin{align}
     m(g_0)-m(g_\infty)=\frac{1}{2(n-1)\omega_{n-1}}\lim_{t\rightarrow\infty}\int R_{g(t)}\ dV_t.\label{constmass}
\end{align}
%Then the mass among all such metrics in the conformal class $(M^n,[g_0])$ is minimized by the unique scalar flat metric in this class.
\end{Cor}

Since for any compact region $K\subset M$, $\lim_{t\rightarrow\infty}\int_K R_{g(t)}\ dV_{g(t)}=0$, we see in \eqref{constmass} that along the Yamabe flow, the difference between the initial mass and the limit mass is accounted for by the total scalar curvature pushed out to infinity by the flow. Such a phenomenon has also been shown by Li to occur for the Ricci flow on asymptotically flat spaces if long-time existence is assumed \cite{YuLi}. We also note that although \cite{CZ} does not show that the mass is constant along the Yamabe flow in general when $n=5$, our convergence results give us additional curvature control that allows us to deduce this as well from their arguments.

\subsection{The main technical estimates}

Theorem \ref{ltegen} is proved straightforwardly from a maximum principle argument using a result of Dilts--Maxwell which states that one can prescribe any strictly negative scalar curvature in the conformal class of AF metrics associated to any asymptotically flat $(M^n,g_0)$ \cite{DM}, along with standard parabolic regularity theory. The bulk of the paper which follows is devoted to the proof of Theorems \ref{uconverge}, \ref{Ypos}, and \ref{Ypos3}. 
%and then Corollary \ref{masscor} follows from these as we described just above.
%and Theorem \ref{Yneg}.

To prove our convergence results in Theorems \ref{uconverge}, \ref{Ypos}, and \ref{Ypos3} we will need to establish the decay of scalar curvature at suitably fast rates. Towards this purpose we again start with bounds on the conformal factor $u(x,t)$ along Yamabe flow using the assumption that there exists a scalar flat, AF metric in the conformal class $(M,[g_0])$ (or equivalently, that $Y(M,[g_0])>0$), which gives us enough control to conclude uniformly bounded scalar curvature for all positive times by local parabolic estimates (see Proposition \ref{longtime}). Then we use the monotonicity of certain integral norms of scalar curvature along with Moser iteration to first show that $\|R\|_{L^\infty}$ tends to zero as $t\rightarrow\infty$.

\begin{Prop}\label{infdec}
Along the Yamabe flow $(M^n,g(t))$ starting from a $C^{k+\alpha}_{-\tau}$ AF manifold with $Y(M^n,[g_0])>0$, $k\geq 3$, and $\tau>0$, we have that
\begin{align}
\sup_{x\in M^n}|R(x,T)|\xrightarrow{T\rightarrow\infty}0.
\end{align}
\end{Prop}

After this we obtain a quantitative decay rate estimate of the $L^\infty$-norm of $R_{g(t)}$ along the flow. This gives rise to the following result, which is an important step towards deriving our desired convergence of the flow.

\begin{Prop}\label{RLinftyDecay}
Let $(M^n,g_0)$ be a $C^{k+\alpha}_{-\tau}$ AF manifold with $Y(M,[g_0])>0$, $k\geq 3$, and $\tau>0$. Then for any $\delta<\frac{\tau}{2}$ there exists $C>0$ such that $\|R\|_{L^\infty}\leq C t^{-1-\delta}$.
\end{Prop}

This decay rate estimate will allow us to derive our first uniform convergence result, Theorem \ref{uconverge}. In fact, it will also be strong enough to allow us to conclude our first weighted convergence result, Theorem \ref{Ypos}, once we further assume on the order of asymptotic flatness that $\tau>1$. However, as mentioned before this condition is too restrictive when $n=3$ because it excludes many manifolds with well-defined mass. So in this case we need another estimate which gives faster decay of the scalar curvature than in Proposition \ref{RLinftyDecay}. By adding the natural conditions of nonnegative, integrable scalar curvature, we can satisfactorily weaken the restriction on $\tau$. Although analogues of the following estimate can be proven for other dimensions as well, we will only state it for dimension $n=3$ since Proposition \ref{RLinftyDecay} already sufficiently covers dimensions $n\geq 4$.

\begin{Prop}\label{n3decay}
In the setting of Proposition \ref{RLinftyDecay}, if $n=3$, $\tau>\frac{1}{2}$, $R_{g_0}\geq 0$, and $R_{g_0}\in L^1(M^n,g_0)$, then for all $\alpha<\frac{3}{2}$ there exists $C>0$ such that $\|R\|_\infty\leq C t^{-\alpha}$.
\end{Prop}

Together, Propositions \ref{RLinftyDecay} and \ref{n3decay} allow us to establish that $u(t)$ must converge to some $u_\infty$ which is asymptotic to $1$ and is a conformal factor corresponding to a scalar flat deformation of $g_0$. We can then conclude that $u_\infty$ must be the conformal factor corresponding to the unique $C^{k+\alpha}_{-\tau}$ scalar flat metric in the conformal class of $g_0$, and proceed to prove the weighted convergence results of Theorems \ref{Ypos} and \ref{Ypos3}. 

\subsection{Organization of the article}

The organization of the article is as follows:
In Section \ref{shorttime} we start by recalling some preliminaries on the short time existence of Yamabe flow on AF manifolds as well as definitions of the relevant weighted H\"older spaces and Sobolev spaces, before proceeding in Section \ref{boundconformal} to prove the bounds on the conformal factors $u(x,t)$ needed in the rest of the paper. These bounds then allow us to prove the general long-time existence result, Theorem \ref{ltegen}, as well as the uniform scalar curvature estimate of Proposition \ref{longtime} when $Y(M^n,[g_0])>0$.
%We will then prove that certain $L^2$ Euclidean-type Sobolev inequalities holds on $(M^n,g(t))$ along the flow. 
%In Section \ref{decaysection}, we will prove the long-time existence. There are two cases that we consider. In the first case when $g_0$ is a general $C^{k+\alpha}_{-\tau}$ AF metric, we will prove Theorem \ref{ltegen}. In the second case when $g_0$ is a $C^{k+\alpha}_{-\tau}$ AF metric satisfying $Y(M^n,[g_0])>0$, we will prove Theorem \ref{longtime}. 
In Section \ref{evoscalar} we study the decay of the scalar curvature under the flow when $Y(M,[g_0])>0$ in order to prove Propositions \ref{RLinftyDecay} and \ref{n3decay} on the decay of scalar curvature. In Section \ref{convsec}, we use these estimates to first prove uniform convergence in Theorem \ref{uconverge}, and then our weighted convergence results of Theorems \ref{Ypos} and \ref{Ypos3}. 
%Finally in Section \ref{massandneg} we show the relation between the mass and our convergence results by proving Corollary \ref{masscor}, and also collect some observations on the long-time behavior of the flow when $Y(M,[g_0])\leq 0$.

\section{Preliminaries}\label{shorttime}

We begin with a brief general discussion of Yamabe flow, before proceeding to introduce background results on the short-time existence of Yamabe flow on asymptotically flat manifolds and conformal deformations.

Suppose that $(M^n,g(t))$ evolves according to the Yamabe flow with initial metric $g_0$, satisfying
\begin{align}
\begin{cases}
\frac{\partial}{\partial t} g=-R g,
\\
g(0)=g_0.
\end{cases}
\end{align}
Since $g(t)$ remains in the same conformal class along the flow, we may write $g(t)=u(t)^{\frac{4}{n-2}}g_0$. Then we have the following relation between $R_g$ and $R_{g_0}$:
\begin{align}
-a(n)\Delta_{g_0} u+R_{g_0} u=R_g u^{\frac{n+2}{n-2}},
\end{align}
where $a(n)=\frac{4(n-1)}{n-2}$. Thus the Yamabe flow can be rewritten as an evolution equation for the conformal factor $u(t)$:
\begin{align}
\frac{\partial}{\partial t} u^{\frac{n+2}{n-2}}=\frac{n+2}{4}\left(a(n)\Delta_{g_0} u-R_{g_0} u\right).\label{Ypde}
\end{align}
Below, we will often denote $N=\frac{n+2}{n-2}$.

\subsection{Analytic preliminaries}\label{analyticprelim}

We first recall some standard function spaces and related definitions used in the analysis and definition of asymptotically flat (AF) manifolds.  See for instance \cite{Bartnik86,DM}.

\begin{Def}\label{funcspaces}
Let $M^n$ be a complete differentiable manifold such that there exists a compact $K\subset M^n$ and a diffeomorphism $\Phi:M^n\backslash K\rightarrow \mb{R}^n\backslash B_{R_0}(0)$, for some $R_0>0$. Let $r\geq 1$ be a smooth function on $M^n$ that agrees under the identification $\Phi$ with the Euclidean radial coordinate $|x|$ in a neighborhood of infinity, and let $\hat{g}$ be a smooth metric on $M^n$ which is equal to the Euclidean metric in a neighborhood of infinity under the identification $\Phi$. Then with all quantities below computed with respect to the metric $\hat{g}$, we have the following function spaces:

The weighted Lebesgue spaces $L^q_\beta(M)$, for $q\geq 1$ and weight $\beta\in\mb{R}$, consist of those locally integrable functions on $M$ such that the following respective norms are finite:
\begin{align}\notag
\|v\|_{L_{\beta}^{q}(M)}=\left\{\begin{array}{ll}{\left(\int_{M}|v|^{q} r^{-\beta q-n} d x\right)^{\frac{1}{q}},} & {q<\infty}, \\ {\operatorname{ess} \sup _{M}\left(r^{-\beta}|v|\right),} & {q=\infty}.\end{array}\right. 
\end{align}

The weighted Sobolev spaces $W^{k,q}_\beta(M)$ are then defined in the usual way with the norms
\begin{align}\notag
\|v\|_{W_{\beta}^{k, q}(M)}=\sum_{j=0}^{k}\left\|D_{x}^{j} v\right\|_{L_{\beta-j}^{q}(M)}.
\end{align}

The weighted $C^k$ spaces $C^k_\beta(M)$ consist of the $C^k$ functions for which the following respective norms are finite:
\begin{align}\notag
\|v\|_{C_{\beta}^{k}(M)}=\sum_{j=0}^{k} \sup _{M} r^{-\beta+j}\left|D_{x}^{j} v\right|.
\end{align}

The weighted H\"{o}lder spaces $C_\beta^{k+\alpha}(M)$, $\alpha\in(0,1)$, consist of those $v\in C_\beta^k(M)$ for which the following respective norms are finite:
\begin{align}\notag
\|v\|_{C_{\beta}^{k+\alpha}(M)}=\|v\|_{C_{\beta}^{k}(M)}+\sup _{x \neq y \in M} \min (r(x), r(y))^{-\beta+k+\alpha} \frac{\left|D_{x}^{k} v(x)-D_{x}^{k} v(y)\right|}{d(x,y)^{\alpha}}.
\end{align}
\end{Def}

\begin{Rem}
The function spaces defined above are independent of the choices of $\hat{g}$ and $r$. In fact, different choices of the metric $\hat{g}$ and the positive function $r$ will produce equivalent norms. Since $\hat{g}$ agrees with the Euclidean metric in a neighborhood of infinity, we will often use $\delta_{ij}$ to denote a choice of metric $\hat{g}$.
\end{Rem}

We can now define our precise notions of asymptotically flat metrics. An asymptotically flat manifold is then a smooth manifold with an asymptotically flat metric.

\begin{Def}[Asymptotically flat metrics]
Given $M^n$ as in Definition \ref{funcspaces}, a metric $g$ is said to be a $W^{k,q}_{-\tau}$ (respectively $C^k_{-\tau}$, $C^{k+\alpha}_{-\tau}$) asymptotically flat (AF) metric if $\tau>0$ and
\begin{align}
g-\hat{g}\in W^{k,q}_{-\tau}(M)\quad\text{(respectively $C^k_{-\tau}(M)$, $C^{k+\alpha}_{-\tau}(M)$)}.
\end{align}
The number $\tau>0$ is called the order of the asymptotically flat metric.
\end{Def}

\subsection{Short-time existence}

We now recall the short-time existence results for Yamabe flow on asymptotically flat manifolds which will be used in this paper. The short-time existence of the Yamabe flow starting from an asymptotically flat manifold and the preservation of asymptotic flatness along the flow have been established by Cheng--Zhu for $C^{2+\alpha}_{-\tau}$ AF metrics (quoted in Theorem \ref{thm:2.1} below). However, we will study certain higher order $C^{k+\alpha}_{-\tau}$ AF metrics along the Yamabe flow; hence we will also check that results analogous to those of Cheng--Zhu hold in these cases as well.

We start by recalling the definition of the particular kind of solution of the Yamabe flow which we will consider throughout. As mentioned earlier, whether the Yamabe flows defined below are unique in general remains open.

\begin{Def}[{\cite[Definition 1.2]{CZ}}]\label{finedef}
We say that $g(t)$ is a fine solution of the Yamabe flow on a complete manifold $(M^n,g_0)$ on a maximal time interval $[0,T_0)$ if $g(t)=u(t)^{\frac{4}{n-2}}g_0$ with $u(0)\equiv 1$ satisfies \eqref{floweq} and for any $T\in(0,T_0)$ there exists $\delta=\delta(T)$ and $C=C(T)$ such that on $[0,T]$, $0<\delta\leq|u(x,t)|\leq C$, $\sup_{[0,T]\times M^n}|\nabla_{g_0} u(x,t)|\leq C$, $\sup_{[0,T]\times M^n}|Rm(g)|(x,t)\leq C$, and moreover either $T_0<\infty$ and $\lim_{t\rightarrow T_0}|Rm|(\cdot,t)=\infty$, or $T_0=\infty$.
\end{Def}

\begin{Rem}
In fact by \cite[Theorem 1]{Ma}, the blowup alternative in the definition above also holds true when rewritten in terms of the scalar curvature: we must have either $T_0<\infty$ and $\lim_{t\rightarrow T_0}|R|(\cdot,t)=\infty$ or $T_0=\infty$. 
\end{Rem}

We now quote Cheng--Zhu's results on the existence of fine solutions to the Yamabe flow starting from $C^{2+\alpha}_{-\tau}$ AF manifolds.

\begin{Thm}[{\cite[Corollary 2.5]{CZ}}] \label{thm:2.1}
If $(M^n,g_0)$ is a $C^{2+\alpha}_{-\tau}$ AF manifold, $\tau>0$, then there exists a fine solution of the Yamabe flow starting from $(M^n,g_0)$ on a maximal time interval $[0,T_0)$ with $T_0>0$.
\end{Thm}

\begin{Thm}[{\cite[Theorem 1.3]{CZ}}]\label{AFpreserve}
Let $u(x,t)$ on $0\leq t<T_0$ be the conformal factor corresponding to a fine solution to the Yamabe flow on a $C^{2+\alpha}_{-\tau}$ AF manifold $(M^n,g_0)$ with $u(\cdot,0)\equiv 1$, and let $v=1-u$. Then $v(x,t)\in C_{-\tau}^{2+\alpha}(M)$. Hence $g_{ij}(t)-\delta_{ij}\in C^{2+\alpha}_{-\tau}(M)$ for $t\in[0,T_0)$, and in particular $(M^n,g(t))$ remains a $C^{2+\alpha}_{-\tau}$ AF manifold along the Yamabe flow for $t\in[0,T_0)$.
\end{Thm}

As mentioned before, we require the analogues of the above two results for AF manifolds with estimates also on higher order derivatives --- in particular, for $C^{k+\alpha}_{-\tau}$ AF manifolds, $k\geq 2$. Clearly if we replace $C^{2+\alpha}_{-\tau}$ with $C^{k+\alpha}_{-\tau}$ in the statement of Theorem \ref{thm:2.1} the statement remains true. Thus we can conclude with the analogue of Theorem \ref{AFpreserve} below. The proof is a straightforward adaptation of Cheng--Zhu's proof of Theorem \ref{AFpreserve} and is presented in Appendix \ref{appendixcka}.

\begin{Thm}\label{AFpreservecka}
Let $u(x,t)$ on $0\leq t<T_0$ be the conformal factor corresponding to a fine solution to the Yamabe flow on a $C^{k+\alpha}_{-\tau}$ AF manifold $(M^n,g_0)$, $k\geq 2$, with $u(\cdot,0)\equiv 1$, and let $v=1-u$. Then $v(x,t)\in C_{-\tau}^{k+\alpha}(M)$. Hence $g_{ij}(t)-\delta_{ij}\in C^{k+\alpha}_{-\tau}(M)$ for $t\in[0,T_0)$, and in particular $(M^n,g(t))$ remains a $C^{k+\alpha}_{-\tau}$ AF manifold along the Yamabe flow for $t\in[0,T_0)$.
\end{Thm}

\section{Bounds on the conformal factor and long-time existence}\label{boundconformal}

In this section we recall the results of \cite{CB,DM} on conformal deformations of asymptotically flat metrics in order to obtain upper and lower bounds on the conformal factor $u(t)$ as it evolves along Yamabe flow. These bounds will then imply the long-time existence of any fine Yamabe flow starting from an asymptotically flat manifold.

\subsection{Conformal deformations of asymptotically flat metrics}\label{conformaldefsection}

Observe that if $\tilde{g_0}=v^{\frac{4}{n-2}} g_0$ is a fixed metric conformal to the initial metric $g_0$ on $M^n$ from which we start a Yamabe flow, then the conformal factor $u(t)$ of \eqref{Ypde} also satisfies for $w=u(t)v^{-1}$,
\begin{align}
\frac{\partial}{\partial t} w^{\frac{n+2}{n-2}}=\frac{n+2}{4}\left(a(n)\Delta_{\tilde{g_0}}w-R_{\tilde{g_0}}w\right),\label{Ypdeconformal}
\end{align}
which is exactly \eqref{Ypde} but with $\tilde{g_0}$ and $w(t)$ replacing $g_0$ and $u(t)$, respectively.

This suggests making an advantageous choice of background metric $\tilde{g_0}$ with which to study \eqref{Ypdeconformal}. First, Dilts--Maxwell showed that for suitable $W^{k,p}_{-\tau}$ AF manifolds, it is always possible to deform to negative scalar curvature \cite{DM}. For our purposes it is more convenient to work with $C^{k+\alpha}_{-\tau}$ AF manifolds, and the analogous statement holds in this setting as well.

%The following result of Dilts--Maxwell is useful for this purpose:
%; it states that given an asymptotically flat manifold, any prescribed strictly negative scalar curvature $R'$ of sufficient decay, namely $R'\in L^p_{-2-\tau}$, can be achieved by some conformal deformation of the original metric:

\begin{Prop}\label{DMthm}
%[c.f. {\cite[Lemma 4.3]{DM}}]
Let $\left(M^{n}, g\right)$ be a $C_{-\tau}^{k+\alpha}$ AF manifold, $k\geq 2$ with $\tau \in(0,n-2)$. Suppose $R' \in C_{-2-\tau}^{k+\alpha}$ satisfies $R'\leq R_g$. Then there exists a positive function $\phi$ with $\phi-1\in C_{-\tau}^{k+\alpha}$ such that the scalar curvature of $g'=\phi^{\frac{4}{n-2}}g$ is $R'$. In particular $g'$ is also a $C^{k+\alpha}_{-\tau}$ AF metric.
\end{Prop}

If instead we want to conformally deform to $R'\equiv 0$, then work of Cantor--Brill \cite{CB} (corrected and completed by Maxwell \cite{Maxwell}) tells us that we can do so for AF metrics belonging to suitable $W^{k,p}_{-\tau}$ classes if and only if they are Yamabe positive \cite{CB,Maxwell}. Again, the analogous statement holds for $C^{k+\alpha}_{-\tau}$ AF manifolds.

\begin{Prop}\label{CBY}
%[c.f. {\cites[Theorem 2.1]{CB}[Proposition 3]{Maxwell}}]
Let $(M^n,g)$ be a $C^{k+\alpha}_{-\tau}$ AF manifold, $k\geq 2$, with
%$p\in\left(1,\frac{2n}{n-2}\right)$
$\tau\in\left(0,n-2\right)$. Then the following are equivalent:
\begin{enumerate}[(1)]
\item We have $Y(M,[g])>0$.
%\begin{align}
%Y(M,[g]):=\inf_{\substack{u\in C^{\infty}_0(M),\\{u\neq 0}}} \frac{4\frac{n-1}{n-2}\int_M|\nabla u|^2+R_g u^2\ dV_g}{\left(\int |u|^{\frac{2n}{n-2}}\ dV_g\right)^{\frac{n-2}{n}}}>0\label{Yamabeconstant}
%\end{align}
\item There exists a positive function $\phi$ with $\phi-1\in C^{k+\alpha}_{-\tau}$ such that $\tilde{g}=\phi^{\frac{4}{n-2}}g$ is conformally equivalent to $\tilde{g}$ and $R_{\tilde{g}}\equiv0$.
\end{enumerate}
\end{Prop}

We will describe how Propositions \ref{DMthm} and \ref{CBY} follow from their $W^{k,p}_{-\tau}$ versions in Appendix \ref{confdefappx}. 

\subsection{Conformal factor bounds}

By making a suitable choice of background metric as detailed earlier, with the results of Section \ref{conformaldefsection} we can obtain some control of the conformal factor $u(t)$ of \eqref{Ypde} along the Yamabe flow of an asymptotically flat metric. In turn we can then achieve some control of the Sobolev constant as defined below. 
%Below, for $C^{k+\alpha}_{-\tau}$ AF manifolds we assume $k\geq 3$.
%It was 2 before but need 3 because when embedding into $C^{k,\alpha}_{-\tau}$ with one derivative less, need twice differentiability of \phi.

\begin{Def}\label{Sobdef}
If $(M^n,g)$ is a $C^0_{-\tau}$ AF manifold, then there exists a smallest constant $C_g>0$ such that for every $u\in W^{1,2}(M,g)$, the following $L^2$ Euclidean-type Sobolev inequality holds:
\begin{align}
\left(\int |u|^{\frac{2n}{n-2}}\ dV_g\right)^{\frac{n-2}{n}}\leq C_g\int|\nabla u|^2\ dV_g.
\end{align}
We call $C_g$ the Sobolev constant of the metric $g$.
\end{Def}

\begin{Rem}
It is well known that $C_{n,e}\leq C_g$, where $C_{n,e}$ is the Sobolev constant of the flat metric on $\mb{R}^n$; see for instance \cite[Proposition 4.2]{Hebey}.
\end{Rem}

We will consider two cases --- first, the general case where $(M^n,g_0)$ is an arbitrary $C^{k+\alpha}_{-\tau}$ AF manifold, and second, the case when moreover $Y(M^n,[g_0])>0$. In the general case, below we show that for any finite time $T>0$, the conformal factor $u(t)$ is bounded away from both $0$ and $\infty$.

\begin{Lem}\label{utimebd}
If $u(x,t)$ is a solution of \eqref{Ypde} corresponding to a fine solution of the Yamabe flow starting from the $C^{k+\alpha}_{-\tau}$ AF manifold $(M^n,g_0)$, $k\geq 2$, %and $\tau\in(0,n-2)$, 
then for any $T$ for which the Yamabe flow exists on $[0,T]$, there exists a $C(T)>0$ depending only on $T$ and $g_0$ such that we have the bounds
\begin{align}
0<C(T)^{-1}\leq u(x,t)\leq C(T)<\infty,
\end{align}
for any $(x,t)\in[0,T]\times M$.
\end{Lem}
\begin{proof}
First, observe that
\begin{align}\notag
\frac{\partial}{\partial t} u=-\frac{n-2}{4}R_{g(t)} u\leq-\frac{n-2}{4}\left(\inf_{x\in M}R_{g(0)}(x)\right) u,
\end{align}
since $\inf_{x\in M} R_{g(t)}(x)$ is nondecreasing under the Yamabe flow. Therefore we have $u(x,t)\leq C(T)<\infty$ if the flow exists on $[0,T]$ for $C(T)>0$ depending only on $T$ and $\inf_{x\in M} R_{g(0)}(x)$.

Next, using Proposition \ref{DMthm} we may write $w=u(t) v^{-1}$ as in \eqref{Ypdeconformal}, with $v$ corresponding to a suitable choice of prescribed $R_{\tilde{g_0}}<0$ with rapid decay. The following maximum principle type argument shall give us the lower bound on $w$, and hence $u$ also, thereby completing the proof. 

Let $U\subset M$ be an open set, and let $U_{t_0}=(0,t_0]\times U$, $\Gamma_{t_0}=(\{0\}\times U)\cup([0,t_0]\times\partial U)$. We claim that for $\epsilon>0$, if the Yamabe flow exists on $[0,t_0]$, then the minimum of $w+\epsilon t$ cannot be achieved on $U_{t_0}$. Otherwise, at such a space-time point $(x,t)$ in $U_{t_0}$, we have $\frac{\partial}{\partial t} (w+\epsilon t)\leq 0$, and
\begin{align*}
    0\geq\frac{\partial}{\partial t}(w+\epsilon t)&=w^{1-N}\left(\frac{n-2}{4}\right)\left(a(n)\Delta_{\tilde{g_0}}(w+\epsilon t)-R_{\tilde{g_0}} w\right)+\epsilon
    \\
    &>w^{1-N}\left(\frac{n-2}{4}\right)a(n)\Delta_{\tilde{g_0}}(w+\epsilon t).
\end{align*}

%\begin{align}\notag
%0&\geq \frac{\partial}{\partial t} (w+\epsilon t)^N=N(w+\epsilon t)^{N-1}\frac{\partial}{\partial t} w+\epsilon N(w+\epsilon t)^{N-1}
%\\
%&>\frac{(w+\epsilon t)^{N-1}}{w^{N-1}}\left(\frac{n+2}{4}\right)\left(a(n)\Delta_{\tilde{g_0}}(w+\epsilon t)-R_{\tilde{g_0}} w\right).\notag
%\\
%&>\frac{n+2}{4}a(n)\frac{(w+\epsilon t)^{N-1}}{w^{N-1}}\Delta_{\tilde{g_0}}(w+\epsilon t).\notag
%\end{align}

But this is impossible, since $\Delta_{\tilde{g_0}}(w+\epsilon t)\geq 0$ at the point we consider, so we have proven the claim.

This implies $w(x,t)\geq\inf_{x\in M} w(x,0)>0$, giving us the desired lower bound. Indeed, suppose at some $(x,t)$ that $w(x,t)<\inf_{x\in M} w(x,0)$; then if we take $\epsilon>0$ sufficiently small we also have $w(x,t)+\epsilon t<\inf_{x\in M} w(x,0)$ at this same space-time point. But $w$ is asymptotic to $1$ at spatial infinity in the interval $[0,t_0]$, so by taking $U$ sufficiently large we see that this means $w+\epsilon t$ achieves a minimum in $U_t$, which we have seen is impossible.
%Should be asymptotic to $1$ uniformly in [0,t_0] hence we can conclude.
\end{proof}

%\begin{Cor}\label{Sobtimecontrol}
%If $u(x,t)$ is the solution of \eqref{Ypde} corresponding to a fine solution of the Yamabe flow starting from the $C^{k+\alpha}_{-\tau}$ AF manifold $(M^n,g_0)$, $k\geq 2$, then for any $T$ for which the Yamabe flow exists on $[0,T]$, there exists a $C(g_0,T)>0$ depending only on $g_0$ and an upper bound for $T$ such that for every $u\in W^{1,2}(M,g(t))$, the following Sobolev inequality holds:
%\begin{align}\notag
%\left(\int |u|^{\frac{2n}{n-2}}\ dV_{g(t)}\right)^{\frac{n-2}{n}}\leq C(g_0,T)\int|\nabla u|^2\ dV_{g(t)}.\label{sobtimecontrolineq}
%\end{align}
%\end{Cor}

Lemma \ref{utimebd} immediately implies the long-time existence of the Yamabe flow under the hypotheses of Theorem \ref{ltegen}.

\begin{proof}[Proof of Theorem \ref{ltegen}]
Suppose on the contrary that $(M^n,g(t))$ is a fine solution of the Yamabe flow starting from a $C^{2+\alpha}_{-\tau}$ AF manifold $(M^n,g_0)$ which exists up to a finite-time singularity $T>0$. But Lemma \ref{utimebd} implies that the conformal factor $u(t)$ is uniformly bounded on $[0,T)$, and by estimating as in the proof of Lemma \ref{65.20} it follows by standard parabolic regularity theory that $u(t)$ is uniformly bounded in $C^{2+\alpha}_0$. Hence the scalar curvature of $(M^n,g(t))$ is uniformly bounded on $[0,T)$, contradicting the blowup alternative for fine solutions of the Yamabe flow given in the remark following Definition \ref{finedef}.
\end{proof}

%We see from Corollary \ref{Sobtimecontrol} that we always have time-dependent control of the Sobolev constant along the Yamabe flow. Below, in the case that $Y(M,[g_0])>0$ we shall see that we in fact have time-independent bounds on the conformal factor $u(x,t)$, and therefore time-independent control of the Sobolev constant.

If we additionally suppose that $Y(M^n,[g_0])>0$, then we can uniformly bound the conformal factor $u(t)$ both from above and from below.

\begin{Lem}\label{uLinftybound}
Let $u(x,t)$ be the solution of \eqref{Ypde} given by Theorems \ref{thm:2.1} and \ref{AFpreservecka} starting from the $C^{k+\alpha}_{-\tau}$ AF manifold $(M^n,g_0)$, $k\geq 2$, %and $\tau\in(0,n-2)$, 
and further suppose that $Y(M^n,[g_0])>0$. Then for any time interval $[0,T]$ on which $u$ is defined and any $(x,t)\in[0,T]\times M$, there exist $C_1,C_2>0$ depending only on $g_0$ such that
\begin{align}
0<C_1\leq u(x,t)\leq C_2<\infty.\label{unibd}
\end{align}
\end{Lem}
\begin{proof}
We already saw in the proof of Lemma \ref{utimebd} how to obtain the lower bound $0<C_1\leq u(x,t)$, so we only need to justify the upper bound.

Write $w=u(t)v^{-1}$ as in \eqref{Ypdeconformal}, with $v$ given by the conformal factor corresponding to $R_{\tilde{g_0}}\equiv 0$ whose existence is given by Proposition \ref{CBY}. Using the same notation as in the proof of Lemma \ref{utimebd},  we claim that for $\epsilon>0$,
%sufficiently small
the maximum of $w-\epsilon t$ cannot be achieved on $U_{t_0}$. Otherwise, at such a point in $U_{t_0}$, we have $\frac{\partial}{\partial t} (w-\epsilon t)\geq 0$, and therefore
\begin{align*}
    0\leq \frac{\partial}{\partial t}(w-\epsilon t)&=w^{1-N}\left(\frac{n-2}{4}\right)a(n)\Delta_{\tilde{g_0}}(w-\epsilon t)-\epsilon
    \\
    &<w^{1-N}\left(\frac{n-2}{4}\right)a(n)\Delta_{\tilde{g_0}}(w-\epsilon t).
\end{align*}

%\begin{align}\notag
%0\leq \frac{\partial}{\partial t} (w-\epsilon t)^N&=N(w-\epsilon t)^{N-1}\frac{\partial}{\partial t} w-\epsilon N(w-\epsilon t)^{N-1}
%\\
%&<Nw^{N-1}\frac{\partial}{\partial t} u\notag
%\\
%&=\Delta (w-\epsilon t).\notag
%\end{align}
But $\Delta (w-\epsilon t)\leq 0$ at such a maximum point, giving a contradiction, thus proving the claim.

This implies the upper bound. For if instead there exists $(x,t)$ such that $w(x,t)>\sup_{y\in M}w(0,y)$, then $u(x,t)-\epsilon t>\sup_{y\in\mb{R}^n}u(0,y)$ for $\epsilon>0$ sufficiently small. Since $u$ is asymptotic to $1$ at spatial infinity, if we take $U$ sufficiently large, then $u-\epsilon t$ achieves a maximum in $U_t$, which we have seen is impossible. 
\end{proof}

Lemma \ref{uLinftybound} then implies the uniform boundednes of the scalar curvature for all positive times by the same argument used to prove Theorem \ref{ltegen}. We will need this fact later.

\begin{Prop}\label{longtime}
Let $(M^n,g_0)$ be a $C^{k+\alpha}_{-\tau}$ AF manifold with $Y(M,[g_0])>0$, $k\geq 3$, and $\tau>0$. Then the Yamabe flow starting from $(M^n,g_0)$ has scalar curvature uniformly bounded in time for all $t>0$.
\end{Prop}

\begin{proof}%[Proof of Proposition \ref{longtime}]
Because of the uniform bounds from Lemma \ref{uLinftybound}, standard parabolic regularity theory applied to the evolution equation of $u(t)$ implies that the scalar curvature $R$ is uniformly bounded for all $t>0$.
\end{proof}

Moreover, as a consequence of Lemma \ref{uLinftybound}, for the Yamabe flow starting from an asymptotically flat manifold with $Y(M^n,[g_0])>0$ we have uniform control of the Sobolev constant of $(M^n,g(t))$ for all positive times. This will allow us to study the convergence of the flow as $t\rightarrow\infty$.

\begin{Cor}\label{Sobunicontrol}
If $u(x,t)$ is the solution of \eqref{Ypde} corresponding to a fine solution of the Yamabe flow starting from the $C^{k+\alpha}_{-\tau}$ AF manifold $(M^n,g_0)$, $k\geq 2$, and moreover $Y(M,[g_0])>0$, then there exists a constant $D=D(g_0)$ such that for any $T$ for which the Yamabe flow exists on $[0,T]$ the following Sobolev inequality holds for every $u\in W^{1,2}(M,g(t))$, :
\begin{align}
\left(\int |u|^{\frac{2n}{n-2}}\ dV_{g(t)}\right)^{\frac{n-2}{n}}\leq D\int|\nabla u|^2\ dV_{g(t)}.\label{sobuniineq}
\end{align}
\end{Cor}

\section{Evolution of the scalar curvature}\label{evoscalar}

We now study the evolution of the scalar curvature $R$ along the Yamabe flow starting from an asymptotically flat manifold in the $Y(M^n,[g_0])>0$ case, in preparation for proving convergence result of Theorem \ref{Ypos}. Recall that at this point we already have the general long-time existence by Theorem \ref{ltegen}, proved in the previous section. Using in this $Y(M^n,[g_0])>0$ setting the uniform $L^2$ Euclidean-type Sobolev inequality \eqref{sobuniineq} for all positive times, we first establish the monotonicity of certain integral norms of $R$, and then proceed to establish the decay rate estimates of Proposition \ref{RLinftyDecay} and \ref{n3decay} on the $L^\infty$ norm of $R$ in time.

%Given a fixed time interval $[0,T]$, we see from Corollary \ref{Sobtimecontrol} that any $C^{k+\alpha}_{-\tau}$ AF manifold $(M^n,g_0)$ with Yamabe flow existing on this time interval must satisfy a uniform $L^2$ Euclidean-type Sobolev inequality \eqref{sobtimecontrolineq} for every time $t\in[0,T]$. Below we will show that this allows us to first control the $L^{\frac{n}{2}}$ norm of the scalar curvature along the Yamabe flow, and then progress to control higher order integral norms of the scalar curvature. Eventually we reach $L^\infty$ control of the scalar curvature, which implies by \cite[Theorem 2.4]{CZ} that the flow can be extended past $[0,T]$, thus proving the long-time existence of the Yamabe flow starting from AF manifolds in general.

\subsection{Monotonicity of some integral norms of the scalar curvature}\label{monoss}

We have the following equations for the evolution of the scalar curvature and the volume form under the Yamabe flow:
\begin{align}\notag
\frac{\partial}{\partial t} R=(n-1)\Delta R+R^2,\quad \frac{\partial}{\partial t} dV_t=-\frac{n}{2}R\ dV_t
\end{align}

Under the hypotheses of Theorem \ref{Ypos}, we can therefore compute the evolution of $\|R\|_{L^p}$ along the Yamabe flow for $p$ sufficiently large.

\begin{Lem}\label{formula}
Let $(M^n,g_0)$ be a $C^{k+\alpha}_{-\tau}$ AF manifold with $Y(M,[g_0])>0$, $k\geq 3$ and $\tau\in(0,n-2)$. Then for all $p>\frac{n}{2+\tau}$, along the Yamabe flow starting from $(M^n,g_0)$ we have
\begin{align}
\frac{d}{dt}\int|R|^p\ dV_t\leq-\frac{4(n-1)(p-1)}{p}\int|\nabla|R|^{\frac{p}{2}}|^2\ dV_t+\left(p-\frac{n}{2}\right)\int|R|^{p}R\ dV_t.\label{eqformula}
\end{align}
%where $C(n,p)=4\frac{(n-1)(p-1)}{p}$.
\end{Lem}
\begin{proof}
This follows from
\begin{align}\notag
\frac{\partial}{\partial t}R^2=(n-1)\Delta R^2-2(n-1)|\nabla R|^2+2R^3,
\end{align}
and 
\begin{align}
\frac{d}{dt}\int |R|^{p}\ dV_t&=\int\frac{p}{2}(R^2)^{\frac{p}{2}-1}\frac{\partial}{\partial t}R^2-\frac{n}{2}|R|^{p}R\ dV_t,\label{timeintegrable}
\end{align}
since by our assumptions on the asymptotic decay of $(M^n,g_0)$ which are preserved along the flow by Theorem \ref{AFpreservecka} (and indeed for all positive times by Theorem \ref{ltegen}) these integral quantities are well-defined.
\end{proof}

Since $C^{k+\alpha}_{-\tilde{\tau}}$ AF manifolds are also $C^{k+\alpha}_{-\tau}$ if $\tilde{\tau}>\tau$, for the purposes of proving Theorem \ref{Ypos} it suffices to assume that $\tau<n-2$, in which case we always have $p>\frac{n}{2+\tau}>1$. Therefore we will always be assuming that $\tau<n-2$ in the rest of this paper. Next we establish that $\|R\|_p$ is a monotonically nonincreasing quantity along the flow, for appropriate $p$.

\begin{Cor}\label{n2mono}
We have that $\frac{d}{dt}\int|R|^{\frac{n}{2}}\ dV_t\leq 0$.
%If $R_{g_0}\geq 0$ then this holds more generally for $p\in\left(\frac{n}{2+\tau},\frac{n}{2}\right]$.
\end{Cor}
%\begin{proof}
%By Lemma \ref{formula} we have, setting $p=\frac{n}{2}$, that
%\begin{align}\notag
%\frac{d}{dt}\int|R|^{\frac{n}{2}}\ dV_t\leq-4\frac{(n-1)(n-2)}{n}\int|\nabla|R|^{\frac{n}{4}}|^2\ dV_t.\label{n2monoeq}
%\end{align}
%The conclusion follows.
%\end{proof}

As a result we also have monotonicity of $\frac{d}{dt}\int|R|^p\ dV_t$, for $p$ close to $\frac{n}{2}$, as well as the integrability in time of certain other $L^p$ norms of $R$ in space, which will be important when we establish the $L^\infty$ decay of $R$.

\begin{Lem}\label{monoeps1}
There exists $\epsilon=\epsilon(g_0)>0$ such that $\frac{d}{dt}\int|R|^p\ dV_t\leq 0$ for all $p\in\left(\frac{n}{2}-\epsilon,\frac{n}{2}+\epsilon\right)$. Moreover for all such $p$,
\begin{align}
    \int_0^\infty \left(\int|R|^{p\frac{n}{n-2}}\ dV_t\right)^{\frac{n-2}{n}}\ dt<\infty.\label{finint}
\end{align}
\end{Lem}
\begin{proof}
Applying the Sobolev inequality \eqref{sobuniineq} to the first term and the H\"{o}lder inequality to the second term on the right-hand side of \eqref{eqformula}, we obtain
\begin{align}
\frac{d}{dt}\int|R|^p\ dV_t&\leq-\frac{C(n,p)}{D}\left(\int|R|^{p\frac{n}{n-2}}\ dV_t\right)^{\frac{n-2}{n}}\label{pident}
\\
&\quad+\left|p-\frac{n}{2}\right|\left(\int|R|^{p\frac{n}{n-2}}\ dV_t\right)^{\frac{n-2}{n}}\left(\int |R|^{\frac{n}{2}}\ dV_t\right)^{\frac{2}{n}}.\notag
\end{align}
Hence there exists an $\epsilon>0$ such that if $p-\frac{n}{2}<\epsilon$ then
\begin{align}
 -\frac{C(n,p)}{K}+\left|p-\frac{n}{2}\right|\|R_{g_0}\|_{L^{\frac{n}{2}}}<0,\label{decineq}
\end{align}
which implies that $\frac{d}{dt}\int|R|^p\ dV_t\leq 0$, since we know that $\frac{d}{dt}\int|R|^p\ dV_t$ is nonincreasing.

For such $p$, we additionally see that \eqref{pident} implies
\begin{align}\notag
    \frac{d}{dt}\int|R|^p\ dV_t+C\left(\int|R|^{p\frac{n}{n-2}}\ dV_t\right)^{\frac{n-2}{n}}\leq 0,
\end{align}
for some $C$ which may depend on $n$, $p$, and $g_0$. Since $\int |R|^p\ dV_t$ is nonincreasing and nonnegative, we may integrate this inequality to deduce \eqref{finint}.
\end{proof}

\begin{Rem}
In order to prove Proposition \ref{RLinftyDecay} we actually need the monotonicity in Corollary \ref{monoeps1} to hold for all $p>\frac{n}{2+\tau}$ (for sufficiently large times). But to justify this fact we will need to demonstrate the decay of $\|R\|_\infty$ first in Proposition \ref{moserprop} below, before returning to this in Corollary \ref{monointall}.
\end{Rem}

In the above discussions, if $p<2$ (such as when $n$ is small) then one may wish to be careful with the $|\nabla|R|^{\frac{p}{2}}|^2$ integrand in \eqref{eqformula}. We check that the discussion in this section leading up to Lemma \ref{monoeps1} still holds in this context in Appendix \ref{lowreg}.

%\begin{Rem}
%Observe that in the above Proposition, $\epsilon>0$ depends on $T$ because the Sobolev constant used in its proof depends on $T$ (or on an upper bound of $T$).  However, if we additionally have $Y(M^n,[g_0])>0$ and therefore the Sobolev inequality \eqref{sobuniineq} with constant independent of time, then we may take $\epsilon>0$ to be a constant depending only $g_0$.
%\end{Rem}

\subsection{Uniform decay of $R$}

Having obtained the decay of appropriate integral norms of $R$, we can now proceed to control the $L^\infty$ norm of $R$ along the flow by Moser iteration arguments similar to those in \cite{DY}. The central ideas are standard, but we include some details for clarity because we need to adapt them to obtain precise control on the decay rate of $R$, similar to the situation for $|\mathrm{Rm}|$ on certain asymptotically flat Ricci flows studied by the first named author in \cite{EC}.

In fact we will need to pass from $L^p$ to $L^\infty$ control of $R$ several times, so we begin with the following estimate with somewhat general assumptions.

\begin{Prop}\label{moserprop}
Let $(M^n,g(t))$ be a Yamabe flow starting from $C^{k+\alpha}_{-\tau}$ AF manifold with $Y(M,[g_0])>0$, $k\geq 3$, and $\tau\in(0,n-2)$, and suppose that 
\begin{alignat}{2}
&\|R\|_{L^q}\leq\alpha_1 t^{-\gamma_1},\quad &&\text{for $q>\frac{n}{2}$},
\\
&\|R\|_{L^{p_0}}\leq\alpha_2 t^{-\gamma_2},\quad &&\text{for $p_0>\frac{n}{2+\tau}$},\label{p_0decay}
\end{alignat}
for some nonnegative constants $\alpha_1,\alpha_2,\gamma_1,\gamma_2$. Then there exists a constant $C=C(n,q,p_0,\alpha_1,\alpha_2,g_0)>0$ such that
\begin{align}
\sup_{x\in M^n}|R(x,T)|\leq C\max\left\{ T^{\frac{1}{p_0}-\gamma_1\frac{q(n+2)}{p_0(2q-n)} -\gamma_2},T^{-\frac{n}{2p_0}-\gamma_2}\right\}.\label{mosereq}
\end{align}
%Moreover, when $Y(M^n,[g_0])>0$ then $C(n,q,\alpha_1,\alpha_2,g_0,T)>0$ above is actually independent of $T$.
\end{Prop}
\begin{proof}
%As mentioned above, the argument is inspired by that in \cite{DY}; a similar argument has been used in \cite[Proposition 4.8]{EC}. 
We selectively denote $f=|R|$ below in order to distinguish the roles that different factors of $|R|$ play. Applying H\"{o}lder's inequality to the second term on the right in \eqref{eqformula}, we see that for any $p>\frac{n}{2+\tau}$, $q>\frac{n}{2}$, and $\delta>0$,
\begin{align}\notag
\frac{1}{p}\frac{d}{dt}\int f^p\ dV_t&\leq-4\frac{(p-1)}{p^2}(n-1)\int|\nabla f^{\frac{p}{2}}|^2\ dV_t
\\
&\quad+\left(\int |R|^q\ dV_t\right)^{\frac{1}{q}}\left(\delta^{-\frac{n}{2q}}\int f^p\ dV_t\right)^{1-\frac{n}{2q}}\notag
\\
&\qquad\left(\delta^{\left(1-\frac{n}{2q}\right)\frac{n}{n-2}}\int f^{\frac{p}{2}\frac{2n}{n-2}}\ dV_t\right)^{\frac{n-2}{n}\frac{n}{2q}}.\notag
\end{align}
Note that the above in fact holds for all $p>1$, assuming the integrability of all terms involved. We restrict to $p>\frac{n}{2+\tau}$ because of the spatial decay of $f=|R|$ on a $C^{k+\alpha}_{-\tau}$ AF manifold, recalling that we assume $\tau<n-2$. Let $\beta=\beta(t)=\alpha_1 t^{-\gamma_1}$, so that $\|R\|_{L^q}\leq \beta(t)$. We apply Young's inequality to the last term on the right to see that
\begin{align}\notag
\frac{1}{p}\frac{d}{dt}\int f^p\ dV_t&\leq-\frac{4(p-1)}{p^2}(n-1)\int|\nabla(f^{\frac{p}{2}})|^2
\\
&\quad+\beta\delta^{-\frac{n}{2q}}\int f^p\ dV_t+\beta\delta^{1-\frac{n}{2q}}\left(\int f^{\frac{p}{2}\frac{2n}{n-2}}\ dV_t\right)^{\frac{n-2}{n}}.\notag
\end{align}
%Let $T>0$ be such that the Yamabe flow is defined on $[0,T]$, and recall that we have the uniform Sobolev inequality $\|u\|_{L^{\frac{2n}{n-2}}}^2\leq D\|\nabla u\|_{L^2}^2$. Here $D>0$ depends on $g_0$ and an upper bound of $T$ in general by Corollary \ref{Sobtimecontrol}, while when $Y(M^n,[g_0])>0$ then $D$ depends only on $g_0$ by Corollary \ref{Sobunicontrol}.  This is the reason for the dependence in time of $C(n,q,\alpha,g_0,T)$ in \eqref{mosereq}; below we will not explicitly track this particular dependence on $T$. 

Applying the Sobolev inequality \eqref{sobuniineq} to the last term on the right above and setting $\delta=\left(\frac{3(p-1)(n-1)}{\beta D p^2}\right)^{\frac{2q}{2q-n}}$, we obtain for all $p\geq p_0$
\begin{align}
\frac{d}{dt}\int f^p\ dV_t+C_{p_0}\int|\nabla(f^{\frac{p}{2}})|^2\ dV_t\leq C_{p,q,\beta}\int f^p\ dV_t,\label{Cstep1}
\end{align}
where $C_{p_0}=\frac{(p_0-1)(n-1)}{p_0}$ and $C_{p,q,\beta}=p\beta\left(\frac{\beta D p^2}{3(p-1)(n-1)}\right)^{\frac{n}{2q-n}}$.  Now define, for $\frac{T}{2}<\tau<\tau'<T$, the function $\psi:[0,T]\rightarrow[0,1]$,
\begin{align}\notag
\psi(t)=
\begin{cases}
0,\quad &0\leq t\leq\tau,
\\
\frac{t-\tau}{\tau'-\tau},\quad &\tau\leq t\leq\tau',
\\
1,\quad&\tau'\leq t\leq T.
\end{cases}
\end{align}
Then we multiply \eqref{Cstep1} by $\psi$ and find that
\begin{align}\notag
\frac{d}{dt}\left(\psi\int f^p\ dV_t\right)+\psi C_{p_0}\int|\nabla(f^{\frac{p}{2}})|^2\ dV_t\leq \left(C_{p,q,\beta}\psi+\psi'\right)\int f^p\ dV_t,
\end{align}
so that integrating and using the fact that $\beta(t)$ is decreasing, for any $\tilde{t}\in[\tau',T]$ we have
\begin{align}
\int f^p\ dV_{\tilde{t}}+C_{p_0}\int_{\tau'}^{\tilde{t}}\int|\nabla(f^{\frac{p}{2}})|^2\ dV_t\ dt\leq \left(C_{p,q,\beta(\tau)}+\frac{1}{\tau'-\tau}\right)\int_\tau^{T}\int f^p\ dV_t\ dt.\label{Cstep2}
\end{align}
We now define for %$p\geq \frac{n}{2}$ and 
$\tau\in[0,T]$,
\begin{align}\notag
H(p,\tau)=\int_\tau^{T}\int f^p\ dV_t\ dt,
\end{align}
and let $\nu =1+\frac{2}{n}$. We claim that for $p>\frac{n}{2+\tau}$ and $0\leq\tau<\tau'\leq T$,
\begin{align}
H(\nu p,\tau')\leq \frac{D}{C_{p_0}}\left(C_{p,q,\beta(\tau)}+\frac{1}{\tau'-\tau}\right)^\nu H(p,\tau)^\nu.\label{Cstep3}
\end{align}
Indeed,
\begin{align}\notag
\int_{\tau'}^{T}\int f^{\nu p}\ dV_t\ dt&\leq\int_{\tau'}^{T}\left(\int f^p\ dV_t\right)^{\frac{2}{n}}\left(\int f^{\frac{p}{2}\frac{2n}{n-2}}\ dV_t\right)^{\frac{n-2}{n}}\ dt
\\
&\leq D \left(\sup_{\tau'\leq t\leq T}\int f^p\ dV_t\right)^{\frac{2}{n}}\int_{\tau'}^{T}\int|\nabla (f^{\frac{p}{2}})|^2\ dV_t\ dt,\notag
\end{align}
so that \eqref{Cstep2} implies the claim. Now we iterate \eqref{Cstep3} to obtain $L^\infty$ control. Let $p_0$ be as assumed in \eqref{p_0decay}, and define
\begin{align}
\eta=\nu^{\frac{2q}{2q-n}},\quad p_k=\nu^k p_0,\quad\tau_k=\frac{T}{2}+(1-\eta^{-k})\frac{T}{2},\quad \Phi_k=H(p_k,\tau_k)^{\frac{1}{p_k}}.\label{iterdef}
\end{align}
We apply \eqref{Cstep3} to see that
\begin{align}
&\Phi_{k+1}=H(\nu p_k,\tau_{k+1})^{\frac{1}{\nu p_k}}\label{moserrec}
\\
&\leq \left(\frac{D}{C_{p_0}}\right)^{\frac{1}{\nu p_k}}\left(C_{p_k,q,\beta(\tau_k)}+\frac{1}{\tau_{k+1}-\tau_k}\right)^{\frac{1}{p_k}}H(p_k,\tau_k)^{\frac{1}{p_k}}\notag
\\
&\leq \left(\frac{D}{C_{p_0}}\right)^{\frac{1}{\nu p_k}}\left(C_1(n,q,p_0)D^{\frac{n}{2q-n}}\beta(\tau_k) ^{\frac{2q}{2q-n}}+C_2(n,q)\frac{\eta}{T}\right)^{\frac{1}{p_k}}\eta^{\frac{k}{p_k}}\Phi_k\notag
\\
&\leq \left(\frac{D}{C_{p_0}}\right)^{\frac{1}{\nu p_k}}\left(C_1(n,q,p_0)D^{\frac{n}{2q-n}}\left(\frac{\alpha_1}{2}\right)^{\frac{2q}{2q-n}}T^{-\gamma_1\frac{2q}{2q-n}}+C_2(n,q)\frac{\eta}{T}\right)^{\frac{1}{p_k}}\eta^{\frac{k}{p_k}}\Phi_k,\notag
\end{align}
where $C_1(n,q,p_0)=p_0^{\frac{2q}{2q-n}}\left(\frac{p_0}{4(p_0-1)(n-1)}\right)^{\frac{n}{2q-n}}$ and $C_2(n,q)=\frac{2}{\eta-1}$. Since $\sum_{k=0}^\infty\frac{1}{p_k}=\frac{1}{p_0}\frac{n+2}{2}<\infty$ and $\sum_{k=0}^\infty\frac{k}{p_k}<\infty$, we can iterate \eqref{moserrec} to obtain
\begin{align}
\sup_{x\in M^n}|R(x,T)|\leq C(n,q,p_0,\alpha_1,D)\max(T^{-\gamma_1\frac{2q}{2q-n}},T^{-1})^{\frac{1}{p_0}\frac{n+2}{2}}\Phi_0.\label{modifystep}
\end{align}
Finally, we have that
\begin{align}\notag
\Phi_0=\left(\int_{\frac{T}{2}}^T\int|R|^{p_0}\ dV_t\ dt\right)^{\frac{1}{p_0}}\leq C(\alpha_2)T^{-\gamma_2+\frac{1}{p_0}}.
\end{align}
Putting things together we obtain the claimed estimate for $|R|$.
\end{proof}

Next, by slightly modifying the proof of Proposition \ref{moserprop} we will establish Proposition \ref{infdec}, which shows that $\|R\|_{L^\infty}$ decays to zero in time. 

%We did not give this modified version of the proof above because we still need to use the exact statement of Proposition \ref{moserprop} several times later on.

\begin{proof}[Proof of Proposition \ref{infdec}] 
%Applying the Sobolev inequality \eqref{sobuniineq} to \eqref{n2monoeq} from the proof of Lemma \ref{n2mono}, we see that
%\begin{align}\notag
%\frac{d}{dt}\int|R|^{\frac{n}{2}}\ dV_t+\frac{C(n)}{D}\left(\int |R|^{\frac{n}{2}\frac{n}{n-2}}\ dV_t\right)^{\frac{n-2}{n}}\leq 0.
%\end{align}
%Integrating in time from $0$ to $\infty$, since $\int|R|^{\frac{n}{2}}\ dV_t$ is monotonic decreasing we find that
Substituting $p=\frac{n}{2}$ into \eqref{finint} from Lemma \ref{monoeps1}, we have
\begin{align}
\int_0^\infty\left(\int |R|^{\frac{n}{2}\frac{n}{n-2}}\ dV_t\right)^{\frac{n-2}{n}}\ dt<\infty.\label{infdeceq}
\end{align}
Moreover, $\int|R|^{\frac{n}{2}\frac{n}{n-2}}\ dV_t$ is uniformly bounded along the flow, since $\|R\|_{L^\infty}$ and $\|R\|_{L^{\frac{n}{2}}}$ are both uniformly bounded by Proposition \ref{longtime} and Corollary \ref{n2mono}, respectively. Hence
\begin{align}
\int_0^\infty\int |R|^{\frac{n}{2}\frac{n}{n-2}}\ dV_t<\infty.\label{holderdecayest}
\end{align}
We now refer back to the proof of Proposition \ref{moserprop}. Choose $q>\frac{n}{2}$ so that $\|R\|_{L^q}\leq\alpha_1$ for some constant $\alpha_1=\alpha_1(g_0)$, which is possible by Lemma \ref{monoeps1}, and set $p_0=\frac{n}{2}\frac{n}{n-2}$. Following the same steps up to \eqref{modifystep} we obtain, with $\gamma_1=0$, that
%to carry out another Moser iteration argument. First we repeat the arguments in that proof up to \eqref{Cstep3}. At this point, unlike in the proof of Proposition \ref{moserprop}, we set $p_0=\frac{n}{2}\frac{n}{n-2}$; in particular note now that $p_0$ may be different from $q$. The remaining definitions in \eqref{iterdef} remain the same. We therefore obtain a slightly different version of \eqref{moserrec}:
%\begin{align}
%\Phi_{k+1}\leq D^{\frac{1}{\nu p_k}}\left(\widetilde{C_1}(n,q)D^{\frac{n}{2q-n}}\left(\frac{\alpha}{2}\right)^{\frac{2q}{2q-n}}T^{-\gamma\frac{2q}{2q-n}}+C_2(n,q)\frac{\eta}{T}\right)^{\frac{1}{p_k}}\eta^{\frac{k}{p_k}}\Phi_k,\label{moserrec2}
%\end{align}
%where $\widetilde{C_1}(n,q)=\left(\frac{1}{2(n-1)}\right)^{\frac{n}{2q-n}}p_0^{\frac{2q}{2q-n}}$ and $C_2(n,q)$ is as before. This time $\sum_{k=0^\infty}\frac{1}{p_k}=\frac{1}{p_0}\frac{n+2}{2}$, so when we iterate \eqref{moserrec2} we find
%\begin{align}\notag
%\sup_{x\in M^n}|R(x,T)|\leq C(n,q,\alpha,D)\max(T^{-\gamma\frac{q}{p_0}\frac{n+2}{2q-n}},T^{-\frac{n+2}{2p_0}})\Phi_0.
%\end{align}
\begin{align}\notag
\sup_{x\in M^n}|R(x,T)|\leq C(n,q,\alpha_1,D)\max(1,T^{-\frac{n+2}{2p_0}})\Phi_0.
\end{align}
%Now we choose $q>\frac{n}{2}$ so that $\|R\|_{L^q}\leq\alpha$ for some constant $\alpha=\alpha(g_0)$, which is possible by Lemma \ref{monoeps1}.
%as in the proof of Theorem \ref{longtime}. 
%So above, $\gamma=0$. 
Then, since
\begin{align}\notag
\Phi_0=\left(\int_{\frac{T}{2}}^T\int|R|^{p_0}\ dV_t\ dt\right)^{\frac{1}{{p_0}}}\xrightarrow{T\rightarrow\infty}0,
\end{align}
we find that indeed $\|R\|_{L^\infty}\rightarrow 0$ in time.
\end{proof}

Since by Lemma \ref{monoeps1} $\|R\|_{L^p}$ is uniformly bounded along the flow for $p$ close to $\frac{n}{2}$, Proposition \ref{infdec} allows us interpolate to obtain decay of integral norms of $R$ as well.

\begin{Cor}\label{prop:4.2}
For $\epsilon=\epsilon(g_0)>0$ as in Lemma \ref{monoeps1} and $p>\frac{n}{2}-\epsilon$, we have $\|R\|_{L^p}\xrightarrow{t\rightarrow\infty} 0$.
\end{Cor}

In particular, we have $\|R\|_{L^{\frac{n}{2}}}\xrightarrow{t\rightarrow \infty} 0$. This allows us to strengthen the range of exponents covered by Lemma \ref{monoeps1}.

\begin{Cor}\label{monointall}
For any $p>\frac{n}{2+\tau}$, we have $\frac{d}{dt}\int|R|^p\ dV_t\leq 0$ for $t\geq T(p)$ sufficiently large. Moreover for all such $p$,
\begin{align}
\left(\int|R|^{p\frac{n}{n-2}}\ dV_t\right)^{\frac{n-2}{n}}=o(t^{-1}).\label{specialdec}
\end{align}
\end{Cor}
\begin{proof}
In the proof of Lemma \ref{monoeps1} observe that for any fixed $p$, the inequality \eqref{decineq} will hold when $\|R_{g_0}\|_{\frac{n}{2}}$ is replaced by $\|R_{g(t)}\|_{\frac{n}{2}}$ for all $t$ sufficiently large, since we now know that $\|R\|_{\frac{n}{2}}$ decays to zero in time. And because we now know that $\int|R|^{p\frac{n}{n-2}}\ dV_t$ is monotonically nonincreasing in time for $t$ sufficiently large, the integrability expressed in \eqref{timeintegrable} implies \eqref{specialdec}.
\end{proof}

We will use these consequences of Proposition \ref{infdec} to strengthen our control of $R$ and obtain our decay rate estimates on $\|R\|_\infty$.

\subsection{Decay rate of $R$}

Although Proposition \ref{infdec} only told us that $\|R\|_\infty$ tends to zero, we can use the consequent improved integral decay estimates of Corollary \ref{prop:4.2} and \ref{monointall} to show that the $L^\infty$ norm of $R$ must in fact decay at a particular rate, thus proving Proposition \ref{RLinftyDecay}.

\begin{proof}[Proof of Proposition \ref{RLinftyDecay}]
By Corollary \ref{monointall}, we have that for any $p>\frac{n}{2+\tau}$ and $\alpha>0$, for sufficiently large times it holds that $\|R\|_{L^{p\frac{n}{n-2}}}\leq\alpha t^{-\frac{1}{p}}$. Choose $p_1,p_2>\frac{n}{2+\tau}$ with $p_1\frac{n}{n-2}>\frac{n}{2}$, and apply the estimate \eqref{mosereq} of Proposition \ref{moserprop} with $q=p_1\frac{n}{n-2}$ and $p_0=p_2\frac{n}{n-2}$. We thus obtain the estimate (for $T$ sufficiently large),
\begin{align}
\sup_{x\in M^n}|R(x,T)|\leq C(n,q,p_0,\alpha_1,\alpha_2,g_0)\max\left\{T^{\frac{n^2+4q}{np_2(n-2q)}},T^{-\frac{n}{2 p_2}}\right\}.\label{intdecayest}
\end{align}
We now compare the two $T$ exponents above to see that we can achieve the decay exponent in the statement of the Proposition. For the second exponent, clearly as $p_2$ approaches $\frac{n}{2+\tau}$ from above the exponent $-\frac{n}{2 p_2}$ approaches $-1-\frac{\tau}{2}$ from above. For the first exponent, notice that for fixed $p_2>0$ we have $\lim_{q\rightarrow\frac{n}{2}^+}\frac{n^2+4q}{np_2(n-2q)}=-\infty$. If $\tau\geq\frac{4}{n-2}$ then $\frac{n}{2}\geq\frac{n}{n-2}\cdot\frac{n}{2+\tau}$, so we will always be able to choose $q>\frac{n}{2}$ sufficiently close to $\frac{n}{2}$ such that $\frac{n^2+4q}{np_2(n-2q)}<-\frac{n}{2 p_2}$. On the other hand, if $0<\tau<\frac{4}{n-2}$ then the following inequality holds at $q_0=\frac{n}{2+\tau}\frac{n}{n-2}>\frac{n}{2}$:
\begin{align*}
\frac{n^2+4q_0}{np_2(n-2q_0)}<-\frac{n}{2p_2}.
\end{align*}
Hence this inequality remains true if we replace $q_0$ with a $q$ slightly larger than $q_0$, finishing the proof.

%Comparing the exponents above, one may check that $-\frac{1}{p}\frac{n}{2}>-\frac{1}{p}\frac{-n^2+2n-4p}{n(n-2p-2)}$ on the interval $\left(\frac{n-2}{2},\frac{n}{2}\right)$, and moreover $\lim_{p\rightarrow\frac{n-2}{2}^+}-\frac{1}{p}\frac{-n^2+2n-4p}{n(n-2p-2)}= -\infty$. Therefore from \eqref{intdecayest} we conclude that
%\begin{align*}
%\sup_{x\in\mathbb{R}^n}|R(x,t)|\leq C(n,q,\alpha,g_0)T^{-\frac{1}{p}\frac{n}{2}},
%\end{align*}
%and since we may take any $p>\frac{n}{2+\tau}$, the result follows.
%When $p=\frac{n}{2}$ we have $-\frac{1}{p}\frac{n}{2}=-\frac{1}{p}\frac{-n^2+2n-4}{n(n-2p-2)}=-1$. And we have
%\begin{align}\notag
%&\frac{d}{dp}\left.\left(-\frac{1}{p}\frac{n}{2}\right)\right|_{p=\frac{n}{2}}>0,
%\\
%&\left.\frac{d}{dp}\left[-\frac{1}{p}\frac{-n^2+2n-4p}{n(n-2p-2)}\right]\right|_{p=\frac{n}{2}}=\frac{n^2+2n-4}{n^2}>0.
%\end{align}
%Thus we can take $p<\frac{n}{2}$ sufficiently close to $\frac{n}{2}$ in \eqref{intdecayest} to obtain the desired estimate on $\|R\|_{L^\infty}$.
\end{proof}

%\begin{Prop}
%We have that $\|R\|_{L^\infty}=o(t^{-1})$.
%\end{Prop}
%\begin{proof}
%By Lemma \ref{specialdeclem}, we have that given any $\alpha>0$, for sufficiently large times it holds that $\|R\|_{L^{\frac{n}{2}\frac{n}{n-2}}}\leq\alpha t^{-\frac{2}{n}}$. Plugging this into the estimate \eqref{mosereq} of Proposition \ref{moserprop} gives the conclusion.
%\end{proof}

As described in the Introduction, the decay in Proposition \ref{RLinftyDecay} is not quite good enough for our purposes when $n=3$, so we need the improved estimate of Proposition \ref{n3decay} for dimension $n=3$ under additional assumptions of nonnegative, integrable scalar curvature. We conclude this section with its proof.

\begin{proof}[Proof of Proposition \ref{n3decay}]
By \cite[Theorem 1.5]{CZ}, under the assumptions of Proposition \ref{n3decay}, the scalar curvature remains integrable along the flow and moreover $\int R\ dV_t$ is nonincreasing; in fact for dimension $n=3$ we have
\begin{align}\notag
    \frac{d}{dt}\int R\ dV_t\leq-\frac{1}{2}\int R^2\ dV_t\leq 0.
\end{align}
So we may interpolate with the $L^\infty$ estimate of Proposition \ref{RLinftyDecay} to find that $\|R\|_p<\infty$ for all $p\geq 1$, and moreover
\begin{align}
    \|R\|_p=O\left(t^{-\left(1+\delta\right)\left(1-\frac{1}{p}\right)}\right),\label{interpest}
\end{align}
for all $\delta<\frac{\tau}{2}$. By our hypothesis on $\tau$ then the evolution equation for $\frac{d}{dt}\int R^p\ dV_t$ is well-defined whenever $p>1$, so we can apply the estimate \eqref{mosereq} of Proposition \ref{moserprop} with $q=p_1\frac{n}{n-2}=3p_1$ and $p_0=p$, for some $p_1>\frac{3}{2+\tau}$ and $p>1$ to be determined. Plugging in $\|R\|_q\leq \alpha_1 t^{-\frac{1}{p_1}}$ and \eqref{interpest} to estimate $\|R\|_{p_0}$, we obtain
\begin{align}\notag
    &\sup_{x\in M^n}|R(x,T)|
    \\
    &\leq C(n,q,p,\alpha_1,1,g_0)\max\left\{T^{\frac{1}{p}-\frac{n(n+2)}{(n-2)p(2q-n)} -\left(1+\delta\right)\left(1-\frac{1}{p}\right)},T^{-\frac{n}{2p}-\left(1+\delta\right)\left(1-\frac{1}{p}\right)}\right\}.
\end{align}
It suffices to show that given any $\alpha<\frac{3}{2}$ we can choose $q,p_0$ so that both exponents of $T$ are less than $-\alpha$. For the first exponent, setting $p=1$ and $n=3$ we obtain
\begin{align}\notag
    1-\frac{15}{2q-n}<-\frac{18}{7}<-\frac{3}{2},
\end{align}
if we set $q=3\cdot\frac{3}{2+\tau}$, using that $\tau>\frac{1}{2}$. For the second exponent, setting $p=1$ gives us exactly $-\frac{3}{2}$. Therefore by setting $p>1$ close to $1$ and $q>3\cdot\frac{3}{2+\tau}$ close to $3\cdot\frac{3}{2+\tau}$ we obtain our desired estimate.
\end{proof}

\section{Convergence of the Yamabe flow}\label{convsec}
%, assuming $g_0$ is a $C^{k+\alpha}_{-\tau}$ metric}

Using the decay rate estimates on the scalar curvature of $(M^n,g(t))$ evolving under the Yamabe flow from an initial Yamabe-positive asymptotically flat metric satisfying the hypotheses of Proposition \ref{RLinftyDecay} or \ref{n3decay}, we proceed in this section to prove Theorem \ref{Ypos}. Below we start with the unweighted convergence of the conformal factor $u(t)$ before proceeding to study its weighted convergence as $t\rightarrow\infty$.

%Throughout, we work with a Yamabe flow $(M^n,g(t))$ starting from a $C^{k+\alpha}_{-\tau}$ AF manifold with $Y(M,[g_0])>0$, $k\geq 3$, and $\tau\in(0,n-2)$ so that the hypotheses of Proposition \ref{RLinftyDecay} hold, and will indicate additional assumptions where needed.

\subsection{Unweighted convergence}

We first show that the decay rate estimate of Proposition \ref{RLinftyDecay} gives us uniform convergence to a limiting continuous function $u_\infty$. 

\begin{Lem}\label{6uLinfty} 
Under a Yamabe flow satisfying the hypotheses of Proposition \ref{RLinftyDecay}, there exists a continuous function $u_\infty(x)>0$ on $M$ such that 
\begin{align}\label{6.70} \|u(x, t)- u_\infty(x) \|_{L^\infty(M) } \leq \frac{C}{t^{\delta}},
\end{align}
and
\begin{align}
u_\infty(x)-1 \rightarrow 0,\quad\text{as }r\rightarrow\infty,\label{6.110} 
\end{align}
for all $\delta<\frac{\tau}{2}$. If the flow additionally satisfies the hypotheses of Proposition \ref{n3decay}, then this holds for all $\delta<\frac{1}{2}$.
(Recall from Definition \ref{funcspaces} that $r\geq 1$ is a smooth function that agrees with $|x|$ in asymptotic coordinates on $M^n$.)
\end{Lem}

\begin{proof}
Since $\frac{\partial}{\partial t} u = -\frac{n-2}{4} R u$,
we have for all $x\in M^n$ that
\begin{align}\notag
u(x, t)= e^{\int_0^t - \frac{n-2}{4}R(x, t) dt   } u(x, 0).
\end{align}
By Proposition \ref{RLinftyDecay}, for any $0<\delta<\frac{\tau}{2}$,
\begin{align}\notag
|R(x, t)| \leq \frac{C}{ (1+ t)^{1+\delta}   }.
\end{align}
Thus for $t_1,t_2>> 1$,
\begin{align}\notag
|u(x, t_1)- u(x, t_2)|\leq & C  |e^{c-\frac{1}{\delta\cdot t_1^{\delta}} }-e^{c-\frac{1}{\delta\cdot t_2^{\delta}} }      | \cdot u(x, 0)\\
\leq & C| \frac{1}{t_2^{\delta}} - \frac{1}{t_1^{\delta}}  |\cdot | u(x, 0)|\rightarrow 0, \ \mbox{as} \ t_1, t_2 \rightarrow \infty.\notag
\end{align}
This implies
 $ u(x, t)$ is a Cauchy sequence in $L^\infty(M)$ as $t\rightarrow \infty.$
We then immediately conclude that there exists a limiting function $u_\infty(x)\in L^\infty(M)$ which satisfies 
\begin{align}\notag \|u(x, t)- u_\infty(x) \|_{L^\infty(M)} \leq \frac{C}{t^{\delta}}.
\end{align}
And since for each time $t>0$, 
\begin{align}\label{6.100} |u(x,t)-1| \leq \frac{C}{r^\tau}, \quad \mbox{as} \ r\rightarrow \infty,\end{align}
we deduce \eqref{6.110}.

If we are in the setting of Proposition \ref{n3decay}, then the entire argument above carries through but with $\delta<\frac{1}{2}$.
\end{proof}

Note that we could not yet conclude above that
\[ |u_\infty(x)-1| \leq \frac{C}{r^\tau}, \quad \mbox{as} \ r\rightarrow \infty,\]
because the bound $C$ in \eqref{6.100} is not uniform in $t>0$. 
But in the next two propositions we will be able to first show that $u_\infty^{\frac{4}{n-2}}g_0$ is a scalar flat metric, and then show that $u_\infty(x)-1$ does indeed satisfy such a spatial decay estimate.

\begin{Prop} \label{65.20}
%Suppose $g_0$ is a %$C^{k+\alpha}_{-\tau}$ AF metric, $k\geq 2$. Then 
Under a Yamabe flow satisfying the hypotheses of Proposition \ref{RLinftyDecay}, for all $0<\alpha'<\alpha$,
%$u(x, t)\rightarrow u_\infty(x)$ in $C^{k+\alpha'}_{loc}$ as $t\rightarrow \infty$.
\begin{align}\label{6.250}
\|u(x, t)- u_\infty(x)\|_{C^{k+\alpha'}_{loc}}\leq \frac{C}{t^{{\delta}}},\end{align}
for any $0<\delta<\tau/2$. Moreover, $u_\infty\in C^{k+\alpha}_{loc}(M)$ and $g_\infty=u_\infty(x)^{\frac{4}{n-2}}g_0$ is a scalar flat metric, i.e.
\begin{align}\notag
 \Delta_{g_0} u_\infty(x)- \frac{1}{a(n)} R_{g_0}(x) u_\infty(x) =0.\end{align}
Recall $a(n)= \frac{4(n-1)}{n-2}$. If the flow additionally satisfies the hypotheses of Proposition \ref{n3decay}, then \eqref{6.250} holds for all $\delta<\frac{1}{2}$.
\end{Prop}

\begin{proof}

First recall that
\begin{align}\notag
\frac{\partial}{\partial t} u= (n-1)   \left( \frac{1}{u^{N-1}}\Delta_{g_0} u -\frac{\frac{1}{a(n)}R_{g_0}}{u^{N-2}}\right). 
\end{align}
Since we have uniform and positive upper and lower bounds on $u$ by Lemma \ref{uLinftybound},
%and the long-time existence result Theorem \ref{ltegen}, 
%\begin{align}\notag
%0<C_1\leq u(x, t)\leq C_2 < \infty, \ \mbox{for all} \ t\in [0, \infty),
%\end{align}
%where $C_1$, $C_2$ depending only on $g_0$.
and since $R_{g_0}\in C^{k-2+\alpha}_0(M) $ we see that
\begin{align}\notag
\left\|\frac{R_{g_0}}{u^{N-2}} \right\|_{L^\infty_{loc}(M)}\leq C,  \ \mbox{for all} \ t\in [0, \infty).
\end{align}
Thus for all $p\geq 1$, on any compact subset $\Omega\subset \subset M$,
\begin{alignat}{2}\notag
 \|\frac{R_{g_0}}{u^{N-2}} \|_{L^{ p}_0(\Omega\times [t_0, t_0+1])} &\leq C, \quad \ &&\mbox{for all} \ t_0\in [0, \infty).
\\
 \|  u\|_{L^\infty_0(\Omega\times [t_0, t_0+1])} &\leq C,\quad  \ &&\mbox{for all} \ t_0\in [0, \infty).\notag
\end{alignat}
We now apply the standard Schauder estimates for parabolic equations to see that $  u$ is in $W^{2,1, p}_{loc}(\Omega\times [t_0, t_0+1]) $ and satisfies 
\begin{align}\notag
&\|  u\|_{W^{2,1, p}_0(\Omega\times [t_0, t_0+1])} \\
\leq& C\left(  \left\|\frac{R_{g_0}}{u^{N-2}} \right\|_{L^{ p}_0(\Omega\times [t_0, t_0+1])}  + \|  u\|_{L^\infty_0(\Omega\times [t_0, t_0+1])} \right) \notag
\\
\leq &C_{2, p},  \ \mbox{for all} \ t_0\in [0, \infty).  \notag
\end{align}
Therefore, by the results on embedding of suitable Sobolev spaces into H\"older spaces on parabolic domains, see e.g. \cite[Theorem 3.4]{ChenYazhe},
\begin{align}\label{6.72}
\|  u\|_{C^{1+\alpha,\frac{1+\alpha}{2}}_0(\Omega\times [t_0, t_0+1])}\leq C_{1,\alpha}.
\end{align}
where we have chosen $p>n+2$ so that $1-\frac{n+2}{p}=\alpha$.
The embedding constant depends on $n, p,\operatorname{diam}(\Omega)^{-1}$, $\hat{g}$ and the length of time interval (which is $1$ in our case). In particular, it is independent of $t_0$.

To derive higher order local regularity of $u$, we use an induction argument to prove that $ u\in C^{k+\alpha,\frac{k+\alpha}{2}}_0(\Omega\times [t_0, t_0+1])  $. 
Assuming $ u\in C^{l+\alpha,\frac{l+\alpha}{2}}_0(\Omega\times [t_0, t_0+1])  $, we want to show that
$ u\in C^{l+2+\alpha,\frac{l+2+\alpha}{2}}_0(\Omega\times [t_0, t_0+1])  $ if $l\leq k-2$.
From \eqref{6.72}, this is already proved for $l=0$ and $1$.
Since $R_{g_0}\in C^{k-2+\alpha}_{loc}$, by a product estimate for parabolic H\"{o}lder spaces we find for such $l$ that
\begin{align}\notag
 \left\|\frac{R_{g_0}}{u^{N-2}} \right\|_{C^{l+\alpha,\frac{l+\alpha}{2}}_0(\Omega\times [t_0, t_0+1])   } \leq C,  \quad \mbox{for all} \ t_0\in [0, \infty).
\end{align}
Therefore by higher order Schauder estimates for parabolic equations,
$  u\in C^{l+2+\alpha,\frac{l+2+\alpha}{2}}_0(\Omega\times [t_0, t_0+1])$ and satisfies

%W^k, p to W^k+2, p: 1/u^{N-2} \in W^{k, p}     holder case: 1/u^{N-2} \in C^{k, \alpha} 

\begin{align}\notag
&\|  u\|_{C^{l+2+\alpha,\frac{l+2+\alpha}{2}}_0(\Omega\times [t_0, t_0+1]) } \\
\leq & C\left(  \left\|\frac{R_{g_0}}{u^{N-2}} \right\|_{C^{l+\alpha,\frac{l+\alpha}{2}}_0(\Omega\times [t_0, t_0+1])}  + \|  u\|_{L^\infty_{0}(\Omega\times [t_0, t_0+1] )} \right)  \notag \\
\leq & C_{l+2, \alpha}, \quad \mbox{for all} \ t_0\in [0, \infty). \notag
\end{align}
The induction argument stops when $l=k-2$.

Hence we have proved that $ u\in C^{k+\alpha,\frac{k+\alpha}{2}}_0(\Omega\times [t_0, t_0+1]))$.
%Thus by embedding $ u\in C^{k-1+\alpha,\frac{k-1+\alpha}{2}}_{0}(\Omega\times [t_0, t_0+1]    )$ for all $k\geq 0$ and $0<\alpha<1$.
Taking $j=0$ and $t=t_0$ in Definition \ref{A1}, we obtain
\begin{align}\notag
\| u  (x, t)\|_{C^{k+\alpha}_{0}(\Omega    ) }  \leq C_{k+\alpha},   \quad \mbox{for all} \ t\in [0, \infty).
 \end{align}
Note that $C_{k, \alpha} $ depends on $\Omega, k, \alpha$, and is independent of $t$,
%Since \eqref{6.27} is valid for any fixed $k\in \mathbb Z^+$ and $\alpha$.
so for every $0<\alpha'<\alpha$ and every sequence $ \{t_j\}\rightarrow \infty$,  there exist a subsequence $ \{t_{j_k}\}\rightarrow \infty$ and a limiting function $\tilde{u}_\infty(x)$ such that 
\begin{align}\label{6.29}
u(x, t_{j_k})\rightarrow \tilde{u}_\infty(x)  \quad \mbox{in} \ C^{k+\alpha'}_{loc},  \quad  \mbox{as} \ t_j\rightarrow \infty,
\end{align}
by the Arzela--Ascoli theorem.
%Rellich compactness theorem. 
But by Lemma \ref{6uLinfty}, we also have
\begin{align}\notag
u(x, t)\rightarrow u_\infty(x) \quad \mbox{in} \ L^\infty(M),  \quad  \mbox{as} \ t\rightarrow \infty.\end{align}
Thus $u_\infty\equiv  \tilde{u}_\infty$,
%It is obvious that $u_\infty(x)\in C^{k+\alpha'}_{loc}$ with its norm bounded by $C_{k,\alpha'}$.
%From above discussion, we have proved that for every sequence $ \{t_j\}\rightarrow \infty$,  there exists a subsequence $ \{t_{j_k}\}\rightarrow \infty$, such that $u(x, t_{j_k})$ converges in $C^{k+\alpha'}_{loc}$ norm to the same limit. 
%By a simple calculus theorem,
and this implies 
\begin{align}\notag%\label{6.25}
u(x, t)\rightarrow u_\infty(x)  \quad \mbox{in} \ C^{k+\alpha'}_{loc},  \quad  \mbox{as} \ t \rightarrow \infty,
\end{align}

%Since \eqref{6.25} is valid for all $k\in \mathbb Z^+$, we have
%(by Sobolev embedding) 
%$u_\infty(x)\in C^\infty_{loc}$ and $\|u(x, t) -u_\infty(x)\|_{C^k_{loc}(M)} \rightarrow 0$ as $t\rightarrow \infty$ for all $k\in \mathbb Z^+$.
%Together with Proposition \ref{RLinftyDecay}, 
Since $k\geq 2$, we can pass the limit in the scalar curvature equation to see that
 $u_\infty(x)^{\frac{4}{n-2}}g_0$ is a scalar flat metric, i.e.
 $$
 \Delta_{g_0} u_\infty- \frac{1}{a(n)} R_{g_0} u_\infty =0.
 $$
 
Finally, we claim that $u_\infty(x)\in C^{k+\alpha}_{loc}$, even though the sequence does not converge in this space as in \eqref{6.29}.
This is because the $C^{k+\alpha}_{0}(\Omega)  $ norm is lower semicontinuous. Namely, for any compact subset $\Omega\subset \subset M$,
\begin{align}\notag
\| u_\infty  (x)\|_{C^{k+\alpha}_{0}(\Omega    ) }  \leq \lim_{t\rightarrow \infty} \| u  (x, t)\|_{C^{k+\alpha}_{0}(\Omega    ) } \leq C_{k,\alpha}.
\end{align}
\end{proof}

We now have enough information to conclude that $u_\infty$ is in fact the conformal factor $\phi$ from Proposition \ref{DMthm} corresponding to the scalar-flat metric $\tilde{g}$.

\begin{Prop} \label{uinftyspatialdecay}
%For all $\tau'< \tau$ and $p\geq 1$, $u_\infty(x)-1 \in W^{k,p}_{-\tau'}(M)$. (??Hope to replace it by $C^{k+\alpha}_{-\tau}(M)$)
We have $u_\infty-1\in C^{k+\alpha}_{-\tau}$. 
In particular,
%as $k\geq 2$ and $p\geq 1$,  
\begin{alignat}{2}
&|u_\infty(x)-1| \leq \frac{C}{r^{\tau}}, \quad &&\mbox{as} \ r\rightarrow \infty,\label{firstderivativeuinfty}
\intertext{and}
&|\partial_j u_\infty(x)| \leq \frac{C}{r^{1+\tau}}, \quad &&\mbox{as} \ r\rightarrow \infty.
\end{alignat}

\end{Prop}

\begin{proof}
By Proposition \ref{CBY} there exists a positive function $\phi$ with $\phi-1\in C^{k+\alpha}_{-\tau}$ such that $\tilde{g}=\phi^{\frac{4}{n-2}}g_0$ scalar flat. Since 
$u_\infty^{\frac{4}{n-2}}g_0$ is also scalar flat,
%metric, i.e.
% $$
% \Delta_0 u_\infty(x)- \frac{1}{a(n)} R_0(x) u_\infty(x) =0.
% $$
%By Theorem B.1 (see also Cantor-Brill \cite{CB}), if $Y(M, g_0)>0$, then there exists a scalar flat $W^{k,p}_{-\tau'}(M)$ AF metric $\tilde{g}$. We can denote $\tilde{g}=v^{\frac{4}{n-2}} g_0$. Here we have used $g_0\in C^{k+\alpha}_{-\tau}\in W^{k,p}_{-\tau'}(M)$ for all $\tau'< \tau$ and $p\geq 1$.
%Then 
%$v-1\in W^{k,p}_{-\tau'}(M)$
the function $w:= u_\infty\cdot \phi^{-1} $ will satisfy 

$$
 \Delta_{\tilde{g}} w(x) =0.
 $$
Moreover, by Lemma \ref{6uLinfty}. $w(x)-1\rightarrow 0$ as $r\rightarrow \infty$.
Thus the maximum principle asserts that $w\equiv 1$, so that $u_\infty -1=\phi-1\in C^{k+\alpha}_{-\tau}$.
\end{proof}

We can now prove Theorem \ref{uconverge}.

\begin{proof}[Proof of Theorem \ref{uconverge}]
Part (\ref{uconverge1}) of Theorem \ref{uconverge} follows immediately from combining the results of Propositions \ref{65.20} and \ref{uinftyspatialdecay}.

To see Part (\ref{uconverge2}), notice that all we have used to prove Part (\ref{uconverge1}) is that $u(t)$ remains uniformly bounded from above and below for all times. So if $Y(M,[g_0])\leq 0$ but $u(t)$ remains uniformly bounded from above and below, then Lemma \ref{6uLinfty} and Proposition \ref{65.20} again show that $u(t)$ converges uniformly to some $u_\infty\in C^{k+\alpha}_{loc}(M)$, asymptotic to $1$ at spatial infinity. But by the conformal invariance of the Yamabe constant, we have
\begin{align}
    Y(M,[g_0])=\inf_{\substack{v\in C^{\infty}_0(M),\\{v\neq 0}}} \frac{4\frac{n-1}{n-2}\int_M|\nabla_{g_\infty} v|^2\ dV_{g_\infty}}{\left(\int |v|^{\frac{2n}{n-2}}\ dV_{g_\infty}\right)^{\frac{n-2}{n}}}=\frac{1}{C_{g_\infty}}>0,
\end{align}
where the $L^2$ Euclidean-type Sobolev constant $C_{g_\infty}$ of $g_\infty$ exists because $u_\infty$ is asymptotic to $1$ at spatial infinity. This contradicts our initial assumption on the Yamabe constant of $(M,[g_0])$. Therefore from the proof of Lemma \ref{utimebd} we see that we must in fact have $\sup_{x\in M}u(x,t)\xrightarrow{t\rightarrow\infty}\infty$.

It remains to show that the $L^2$ Euclidean-type Sobolev constant of $(M,g(t))$ also blows up as $t\rightarrow\infty$. If not, then we have uniform control of this constant, which allows us to carry out nearly all the arguments of Section \ref{evoscalar} (without any assumption on the boundedness of $u(t)$) until Proposition \ref{infdec}, where we passed from \eqref{infdeceq} to \eqref{holderdecayest} using the $L^\infty$ bound of $R$, which was a consequence of the uniform bounds on $u$ established when $Y(M,[g])>0$. But we can recover the boundedness of $R$ simply with the bound on the Sobolev constant by plugging the uniform bound of $\|R\|_{L^p}$ for $p$ slightly larger than $\frac{n}{2}$ from Lemma \ref{monoeps1} into the estimate of Proposition \ref{moserprop}. In fact the constant in that estimate only depends on the Sobolev constant bound $D$ of $g_0$, 
%in the sense that the Sobolev constant bound $D$ depended on $g_0$, 
%In fact the constant in that estimate only depends on $g_0$ in the sense that the Sobolev constant bound $D$ depended on $g_0$, 
so since we are now assuming control of the Sobolev constant we obtain the same estimates for $\|R\|_\infty$ on the intervals $[T,T+1]$ for all $T\in\mathbb{Z}_{\geq 0}$, and hence uniform bound of $R$ for all times. As a result, we can continue through the same arguments we used to prove $u(t)$ is unbounded earlier. Note that although Proposition \ref{65.20} uses the boundedness of $u(t)$, we do actually obtain this beforehand when we integrate in proving Lemma \ref{6uLinfty}.
\end{proof}

\subsection{Weighted convergence}

We will now prove the weighted convergence of Yamabe flows starting from AF manifolds with $Y(M^n,[g_0])>0$, and need to assume for the rest of this section that $\tau>1$ as in Theorem \ref{Ypos}, or that $n=3$, $\tau>\frac{1}{2}$, $R_{g_0}\geq 0$, and $R_{g_0}\in L^1(M^n,g_0)$ as in Theorem \ref{Ypos3}. The difference between the proofs of Theorems \ref{Ypos} and \ref{Ypos3} occurs with the distinct decay rate estimates, Propositions \ref{RLinftyDecay} and \ref{n3decay}, required to prove Proposition \ref{Rmixdecay} below in the respective settings of the two Theorems. In particular, the decay rate of $R$ from Proposition \ref{RLinftyDecay} alone is not enough to give us the estimates of Lemma \ref{esth} when $n=3$ and $\tau>\frac{1}{2}$, so for this case we must bring in the improved estimate of Proposition \ref{n3decay}.

Using a strategy similar to that used to prove \cite[Lemma 5.2]{YuLi} in a Ricci flow setting, we now prove below Proposition \ref{Rmixdecay}, a decay estimate for the scalar curvature in both space and time. In fact, once we have this, the rest of the proofs of Theorems \ref{Ypos} and \ref{Ypos3} are the same, and we will complete the proofs of both together at the end of this section.

\begin{Prop}\label{Rmixdecay} Under the hypotheses of either Theorem \ref{Ypos} or Theorem \ref{Ypos3}, for any $\tau'<\tau$ there exists some $\delta_0>0$ and $C>0$ depending only on $g_0$ such that 
\begin{align}
|R(x, t)| \leq \frac{C}{r^{\tau'} (1+ t)^{1+\delta_0}   }
\end{align}
for all $(x, t)\in M\times [0, \infty)$ along the flow.
\end{Prop}

%\begin{Rem} Later, when we wish to study the ADM mass along the flow, as was done in the convergence result of \cite[Theorem 5.1]{YuLi} for a Ricci flow setting and in \cite[Theorem 1.3]{CZ} for short-time Yamabe flow, we will need to restrict to $\tau'>\frac{n-2}{2}$ so that the mass is well-defined, by \cite{Bartnik86}. This is only weaker than the assumption $\tau>1$ in Proposition \ref{Rmixdecay} when $n\geq 4$, which as mentioned earlier is why this additional condition appears in the statement of Corollary \ref{masscor}
%in Proposition \ref{Rmixdecay}, because we will not consider the monotonicity of the mass on AF manifolds. It is well-known that by \cite{Bartnik86} $\tau'>\frac{n-2}{2}$ ensures the mass is well-defined on an AF manifold.)
%\end{Rem}

\begin{proof}
Given $\tau'<\tau$, choose $\sigma_0,\sigma_1$ such that $\tau'<\sigma_1<\sigma_0<\tau$. Let $\delta=\frac{\sigma_0}{2}$, and 
%Let $\sigma_0:= 2\delta$, where $\delta$ is any number strictly smaller than $\frac{\tau}{2}$ from Proposition \ref{RLinftyDecay}. Thus $\sigma_0$ is any number strictly smaller than $\tau$. Choose $\sigma_1$ s.t. $\tau'<\sigma_1<\sigma_0<\tau$.
consider $D:= \{ (x, t)\in M\times [0, \infty): r(x)\geq t^{a}    \}$ where $a>1/2$ is to be determined.

For $(x, t)\neq D$, by Proposition \ref{RLinftyDecay}
\begin{align}\notag
|R|\leq C t^{-1-\delta}\leq   C t^{-1-\eta} r^{- \tau'}  
\end{align}
for some $\eta>0$, when $a >1/2$ but is close enough to $1/2$.
%for some $\delta_0>0$ when $a >1/2$ but is close enough to $1/2$.

We now claim the estimate 
\begin{align}\notag
|R| \leq  C  r^{-2- \sigma_1}  \quad \mbox{for} \ (x,t)\in D
\end{align}
 holds. 
In fact, we define 
$h:= r^{4+2\sigma_1 }$ and $w= h \cdot |R|^2$. Then $w$ satisfies the evolution equation
 \begin{align}\label{6.4}
 \left(\frac{\partial}{\partial t} -(n-1)\Delta \right)w\leq(n-1)(B w- 2\nabla \log h \cdot \nabla w)+ 2 h |R|^3.
\end{align}
Here $B:= \frac{ (2|\nabla h|^2- h \Delta h)    }{h^2  }$.
From Lemma \ref{esth} below, for any $\beta<\frac{1}{2}$ we can bound $B$ uniformly by $C(r^{-2}+ r^{-1} t^{-\beta}) \leq C t^{-\frac{a}{2}-\frac{3}{4}}$ on $D$, by choosing $\beta$ suffiently close to $\frac{1}{2}$.
By Proposition \ref{RLinftyDecay}, 
\begin{align}\notag
h|R|^3\leq w\cdot   t^{-1-\delta}= w\cdot  t^{-1-\frac{\sigma_0}{2}}.
\end{align}
Thus \eqref{6.4} becomes
 \begin{align}\notag
 \left(\frac{\partial}{\partial t} -(n-1)\Delta \right)w\leq  - 2(n-1)\nabla \log h \cdot \nabla w+ C w\cdot  t^{-1-\delta'_0}.
\end{align}
where $\delta'_0 := \min\{\frac{a}{2}-\frac{1}{4},  \frac{\sigma_0}{2}  \}>0   $.

On $\partial D$, we have
\begin{align}\notag
|R|\leq C t^{-1-\frac{\sigma_0}{2}}= Cr^{   -(1+\frac{\sigma_0}{2})/a}\leq C r^{-2-\sigma_1}
\end{align}
for  $a$ sufficiently close to $1/2$.
So we can apply the maximum principle on the noncompact manifold as in \cite[Theorem 12.14]{ChowChu} to conclude the claim holds. (See also the statement in \cite[Theorem 2.1]{YuLi}.)

Therefore on $D$ we have
\begin{align}\notag
|R|\leq C r^{-2-\sigma_1} \leq C t^{-1-\delta_0} r^{-\tau'},
\end{align}
for some $\delta_0>0$ and $a$ sufficiently close to $1/2$.
\end{proof}

\begin{Lem}\label{esth} Under the hypotheses of either Theorem \ref{Ypos} or Theorem \ref{Ypos3}, for $r(x)>>1, t>>1$ and any $\beta<\frac{1}{2}$,
\begin{align}
\label{h1}|\nabla_{g(t)} h|_{g(t)}\leq &C |\nabla_{g_0}h|_{g(0)} \leq \frac{Ch}{r}, \\
\label{h2}| \Delta_{g(t)} h |\leq & \frac{Ch}{r} \cdot \left(\frac{1}{r}+\frac{1}{ t^\beta }\right).
\end{align}

\end{Lem}

\begin{proof}
\eqref{h1} is straightforward as  $u$ is uniformly bounded by Lemma \ref{uLinftybound} and
\begin{align}\notag
|\nabla_{g(t)} h|_{g(t)}= u^{-\frac{2}{n-2}} |\nabla_{g_0}h|_{g(0)}.
\end{align}
%(Is this bound $C(T)$ in Lemma 2.10 finite as $T\rightarrow \infty$? The uniform bound is needed here. Can we assume fine solution throughout this paper?)
Regarding \eqref{h2}, we consider
\begin{align}\label{scalareqnh}
%\Delta_{g(t)} h = u^{-\frac{4}{n-2}}\Delta_{g_0} h + \frac{2n}{n-2} u^{\frac{n+2}{n-2}} g_0^{ij} \partial_i h \partial_j u.
\Delta_{g(t)} h = u^{-\frac{4}{n-2}}\Delta_{g_0} h + 2 u^{-\frac{n+2}{n-2}} g_0^{ij} \partial_i h \partial_j u.
\end{align}
We know from calculation
\begin{align}\label{hDelta} 
\Delta_{g_0} h\leq \frac{Ch}{r^2}, 
\end{align}
and
\begin{align}\label{hgrad}
\partial_{i} h\leq \frac{Ch}{r}. 
\end{align}
The next step is to estimate $\partial_j u$.
From the proof of Proposition \ref{65.20}, formula \eqref{6.250}, we see
 \begin{align}\notag
\|u(x, t)- u_\infty(x)\|_{C^{k+\alpha'}_{loc}}\leq \frac{C}{t^{{\delta}}},\end{align}
for any $0<\delta<\tau/2$.
On the other hand, by \eqref{firstderivativeuinfty} from Proposition \ref{uinftyspatialdecay}, we have 
$|\partial_j u_\infty(x)| \leq \frac{C}{r^{1+\tau}}$. Recall that in the $n=3$ case where we imposed additional assumptions on $R_{g_0}$ we can replace the restriction on $\delta$ by $0<\delta<\frac{1}{2}$, and this also holds when $n\geq 4$ since we assumed in this case that $\tau>1$.

So putting these estimates together we can obtain for any dimension $n\geq 3$ and any $\beta<\frac{1}{2}$ that
 \begin{align}\label{nablau}
|\partial_j u(x, t)|\leq  C \left(\frac{1}{t^{\beta}}+ \frac{1}{r^{1+\tau}}\right)\leq C \left(\frac{1}{t^{\beta}}+ \frac{1}{r}\right).
  \end{align}
%for any $0<\delta<\tau/2$.
%Since $\tau>\frac{n-2}{2}$, we can take $\delta>1/2$ when $n\geq 4$. 
 
%\begin{align}\notag \label{nablau}
%|\partial_j u(x, t)|\leq   C (\frac{1}{t^{1/2}}+ \frac{1}{|x|}).
%   \end{align}
% This gives us the estimate of $\partial_j u$. 
Then, plugging \eqref{hDelta}, \eqref{hgrad}, \eqref{nablau} into \eqref{scalareqnh}, we complete the proof of \eqref{h2}.

%Note the original statement of \cite{GilbargTrudinger}[Theorem 7.27] assumes $u\in W^{2,p}_0(\Omega)$, but the bounds in the inequality is independent of $\Omega$ and thus can be easily generalized to $W^{2,p}_0(\mathbb R^n)$. Choose $\epsilon = t^{\frac{1+ \delta_0 }{2} }$. Then  \begin{align}\notag \|\nabla u(x,t) \|_{L^{p}(B(y, 1))}\leq C  t^{-\frac{1+ \delta_0 }{2} }. \end{align}

\end{proof}

Now that we have weighted decay of the scalar curvature, by integrating we can start to prove the weighted convergence of the conformal factor $u(t)$.

\begin{Prop}\label{6uC0tau} Under the hypotheses of either Theorem \ref{Ypos} or Theorem \ref{Ypos3}, for any $\tau'<\tau$ 
%there exists a limiting function $u_\infty(x)$ with $u_\infty(x)-1\in C^0_{-\tau'}$ s.t. 
we have
$$ u(x, t)- u_\infty(x) \rightarrow 0 \quad \mbox{in} \ C^{0}_{-\tau'},  \quad  \mbox{as} \ t\rightarrow \infty,
$$
where $u_\infty\in C^{k+\alpha}_{-\tau}$ is as in Proposition \ref{uinftyspatialdecay}. Moreover, there is a $\delta_0>0$ such that
\begin{align}\label{6.7} \|u(x, t)- u_\infty(x) \|_{C^{0}_{-\tau'}} \leq \frac{C}{t^{\delta_0}}.
\end{align}

\end{Prop}

\begin{proof}
Since $\frac{\partial}{\partial t} u = -\frac{n-2}{4} R u$,
we have 
\begin{align}\label{6.12} 
u(x, t)= e^{\int_0^t - \frac{n-2}{4}R(x, t) dt   } u(x, 0).
\end{align}
By Proposition \ref{Rmixdecay},
\begin{align}\label{6.13} 
|R(x, t)| \leq \frac{C}{r^{\tau'} (1+ t)^{1+\delta_0}   }.
\end{align}
Thus for $t_1,t_2>> 1$,
\begin{align}\notag
|u(x, t_1)- u(x, t_2)|\leq & C e^{\frac{1}{r^{\tau'}}} \cdot |e^{c-\frac{1}{\delta_0 \cdot t_1^{\delta_0}} }-e^{c-\frac{1}{\delta_0 \cdot t_2^{\delta_0}} }      | \cdot u(x, 0)\\
\leq & C\frac{1}{r^{\tau'}} \cdot | \frac{1}{t_2^{\delta_0}} - \frac{1}{t_1^{\delta_0}}  |\cdot | u(x, 0)|.\notag
\end{align}
This implies

\begin{align}\notag
r^{\tau'} \cdot |u(x, t_1)- u(x, t_2)|\leq C   | \frac{1}{t_2^{\delta_0}} - \frac{1}{t_1^{\delta_0}}|\rightarrow 0, \ \mbox{as} \ t_1, t_2 \rightarrow \infty.
\end{align}
Note that $ u(x, t)-1\in C^0_{-\tau'} $ for every $t>0$ by \cite[Theorem 1.3]{CZ}. Thus $ u(x, t)-1 $ is a Cauchy sequence in $C^0_{-\tau'}$ as $t\rightarrow \infty.$
We may then conclude that
%there exists a limiting function $u_\infty(x)-1\in C^0_{-\tau'}$, and it satisfies 
\begin{align}\notag \|u(x, t)- u_\infty(x) \|_{C^{0}_{-\tau'}} \leq \frac{C}{t^{\delta_0}}.
\end{align}
\end{proof}

Before proceeding to the proof of Theorem \ref{Ypos}, we need to establish some uniform control of the conformal factor $u(t)$ in some parabolic weighted Sobolev and H\"{o}lder spaces.

\begin{Lem}\label{6lemv} 
%Suppose $g_0$ is a $C^{k+\alpha}_{-\tau}$ AF metric, $k\geq 2$,  $n\geq 4$. Fix $\tau'$ as in Proposition \ref{Rmixdecay}.  \\
Under the hypotheses of either Theorem \ref{Ypos} or Theorem \ref{Ypos3}, and for any $\tau'<\tau$:
\begin{enumerate}[(1)]
\item\label{6lemv1} For any $p>1$, there exists a constant $C>0$ such that for any $0<\tau_1<\tau'$ and all $t_0\geq 0$.
   \begin{align}\notag
\| u-1\| _{W^{k,k/2, p}_{-\tau_1 }(M\times [t_0, t_0+1 ]  )}\leq C. \end{align}
\item\label{6lemv2} There exists a constant $C>0$, such that for all $t_0\geq  0$, 
\begin{align}\notag
\| u-1\| _{C^{k+\alpha, \frac{k+\alpha}{2}}_{-\tau' }(M\times [t_0, t_0+1 ]  )}\leq C. \end{align} \\
Here $C$ is independent of time $t_0$.
\end{enumerate}
\end{Lem}
\begin{proof}
%By Theorem 5.3 in \cite{CZ} \begin{align}\notag \| u-1\| _{W^{2, 1, p}_{-\tau }(M\times [t_0, t_0+1 ]  )}\leq  C\left(  \| R_0 \|_{L^p_{-\tau-2}(M \times [t_0, t_0+1 ] )}   + \|  u-1\|_{L^p_{-\tau}(M\times [t_0, t_0+1 ]  )} \right). \end{align} 

 To estimate $ \| u-1\| _{W^{k,k/2, p}_{-\tau_1 }(M\times [t_0, t_0+1 ]  )}$, we observe that 
 $u-1$ satisfies the equation 
 \begin{align}\notag
 N u^{N-1} \frac{\partial}{\partial t} (u-1) =\Delta_{g_0} (u-1) - a(n) R_{g_0} (u-1) + a(n) R_{g_0},
\end{align}
where $a(n)= \frac{n-2}{4(n-1)}$. By Proposition \ref{A4} in the Appendix, 
\begin{align}\notag
&\| u-1\| _{W^{k,k/2, p}_{-\tau_1 }(M\times [t_0, t_0+1 ]  )}\\
\leq & C\left(  \| R_{g_0} \|_{W^{k-2,k/2-1, p}_{-\tau_1-2}(M\times [t_0, t_0+1 ]   )}   + \|  u-1\|_{L^p_{-\tau_1}(M\times [t_0, t_0+1 ]  )} \right).\notag
 \end{align} 
Since $R_{g_0}\in C^{k-2+\alpha}_{-\tau-2}(M)$ and $R_{g_0}$ is independent of $t$, 
   \begin{align}\notag
\| R_{g_0} \|_{W^{k-2, k/2-1, p}_{-\tau_1-2}(M \times [t_0, t_0+1 ] )} \leq C.
 \end{align} 
Thus it suffices to prove that $ \|  u-1\|_{L^p_{-\tau_1}(M\times [t_0, t_0+1 ]  )}\leq C$.

Recall \eqref{6.12} and \eqref{6.13} from the proof of Theorem \ref{6uC0tau}:
\begin{align}\notag
u(x, t)= e^{\int_0^t -\frac{n-2}{4} R(x, t) dt   } u(x, 0),
\end{align}
and
\begin{align}\label{t*}
|R(x, t)| \leq \frac{C}{r^{\tau'} (1+ t)^{1+\delta_0}   }.
\end{align}
We also know from the AF assumption on $g_0$ that
\begin{align}\notag
|u(x, 0)- 1|\leq  \frac{C}{r^\tau}.
\end{align}
Thus 
\begin{align}\label{6.73} 
|u(x, t)- 1|  \leq&  e^{\int_0^t -\frac{n-2}{4} R(x, t) dt   } \cdot |u(x, 0)-1|+ |e^{\int_0^t - \frac{n-2}{4}R(x, t) dt   } -1|  \\
\leq & C  \cdot |u(x, 0)-1| + C \left|\int_0^t - R(x, t) dt    \right|  \notag \\
\leq & \frac{C}{r^{\tau'}}.\notag
\end{align}
Above in the second line of \eqref{6.73}, we have used $|e^s-1|\leq e^{|\bar{s}|}\cdot |s|$, where $\bar{s}$ is some point in $[0, s]$, so that $|e^s-1|\leq e^{|s|} \cdot |s|$. Letting $s= \int_0^t - \frac{n-2}{4}R(x, t) \ dt  $ so that $|s|\leq \frac{-C}{ r^{\tau'} (1+ t')^{\delta_0}   }\big|_{t'=0}^{t}\leq C$ by \eqref{t*} gives us the estimate.
%,  As a result, 
%\[ |e^{\int_0^t - \frac{n-2}{4}R(x, t) dt   } -1|\leq C |\int_0^t - R(x, t) dt    | .\]\\

%Thus by the definition of weighted $L^p$ spaces
%   \begin{align}\notag
% \|  u-1\|_{L^p_{-\tau_1}(M\times [t_0, t_0+1 ]  )}\leq C, \end{align} 
% for any $\tau_1<\tau'$. 
As a result we see that $u-1\in C^0_{-\tau'}(M\times [t_0, t_0+1 ]  )\hookrightarrow L^p_{-\tau_1}(M\times [t_0, t_0+1 ]$, for any $\tau_1<\tau'$,  completing the proof of (\ref{6lemv1}).

To prove (\ref{6lemv1}), again by Proposition \ref{A4} in the Appendix we have
%by \eqref{6.73}, $\|u- 1\|_{C^0_{-\tau'}(M\times [t_0, t_0+1 ]    )   } \leq C$.
%Using $R_0\in C^{k-2+\alpha}_{-\tau-2}(M)$ and
\begin{align}\notag
&\| u-1\| _{C^{k+\alpha, \frac{k+\alpha}{2}}_{-\tau' }(M\times [t_0, t_0+1 ]  )}  \\
\leq & C\left(  \| R_{g_0} \|_{C^{k-2+\alpha, \frac{k-2+\alpha}{2}}_{-\tau'-2}(M\times [t_0, t_0+1 ]   )}   + \|  u-1\|_{C^0_{-\tau'}(M\times [t_0, t_0+1 ]  )} \right).\notag
 \end{align}
 Since $u- 1\in C^0_{-\tau'}(M\times [t_0, t_0+1 ]    )$ and $R_{g_0}\in C^{k-2+\alpha}_{-\tau-2}(M)$ we obtain our desired estimate.
%we complete the proof of (2).
\end{proof}

Now we are ready to prove Theorems \ref{Ypos} and \ref{Ypos3}. 

%First we notice that the statement of Theorem \ref{Ypos} is equivalent to the following theorem.
%\begin{Thm}\label{thm:6.8'} Suppose $g_0$ is a $C^{k+\alpha}_{-\tau}$ AF metric, $k\geq 2$, $n\geq 4$ and
%$\frac{n-2}{2}<\tau$. Then for any $\tau'<\min\{\tau,n-2\}$,
% $\frac{n-2}{2}<\tau<n-2$. Then for any $\frac{n-2}{2}<\tau'<\tau$\\

%$$u(x, t)  -1\rightarrow u_\infty(x)-1  \quad \mbox{in} \ C^{k+\alpha}_{-\tau'}, \quad  \mbox{as} \ t\rightarrow \infty.$$ 
%More precisely, 
%for each $k\in \mathbb Z^+$,
%there exists constants $C_{k,\alpha}>0$ and $\delta_0>0$ s.t.
%\begin{align}\notag\|u(x, t) - u_\infty(x)\|_{ C^{k+\alpha}_{-\tau'} }\leq %\frac{C_{k,\alpha}}{t^{ \delta_0  }},  \quad  \mbox{as} \ t\rightarrow \infty.
%\end{align}
%\end{Thm}

\begin{proof}[Proof of Theorems \ref{Ypos} and \ref{Ypos3}]
%of Theorem \ref{thm:6.8'}. 
We will bootstrap to inductively prove higher global regularities of $\phi$ in weighted parabolic H\"{o}lder spaces.
More precisely, we want to show the following estimate: for all $l\leq k$ and $\tau'<\tau$ (assuming as throughout this section that $\tau<n-2$),
\begin{align}\label{6.61} \|   \phi(x, t)\|_{C^{l+\alpha,\frac{l+\alpha}{2}}_{-\tau'}(M\times [t_0, t_0+1])  } \leq  \frac{C}{t_0^{\delta_0}}. \end{align}

\underline{\bf{Claim 1:}} \eqref{6.61} holds for $l=0$.

Let $\phi(x, t)= u(x, t)-u_\infty(x)  $.
Subtracting the two equations that $u(x, t)$ and $u_\infty(x)$ satisfy,
\begin{align}\notag
\frac{\partial}{\partial t} u^N= (n-1)N   [  \Delta_{g_0} u -  \frac{1}{a(n)}R_{g_0}(x) u] 
\intertext{and}
\Delta_{g_0} u_\infty(x)- \frac{1}{a(n)} R_{g_0}(x) u_\infty(x) =0,\notag
\end{align}
we obtain
\begin{align}\displaystyle
\frac{\partial}{\partial t}\phi(x, t) 
 = \frac{n-1}{u^{N-1}} [ \Delta_{g_0} \phi(x, t) -  \frac{1}{a(n)} R_{g_0}(x) \phi(x, t) ].\label{phieqn}
\end{align}
By \eqref{6.7} of Proposition \ref{6uC0tau}, there is a $\delta_0>0$ such that
\begin{align}%\label{6.23}
\|  \phi(x, t) \|_{C^0_{-\tau'}(M )} \leq \frac{C}{t^{\delta_0}}.\notag\end{align}
Therefore for all $  p\geq 1    $ and $\tau_1<\tau'$,
\begin{align}\label{6.21} \|  \phi(x, t) \|_{L^p_{-\tau_1}(M\times [t_0, t_0+1 ]  )} \leq \frac{C}{t_0^{\delta_0}},\end{align}
where $C$ depends only on $p\geq 1$ and $g_0$ (but is independent of $t_0$).

Since $u(x, t)-1 \in C^0_{-\tau'}(M)$ and $0<C_1\leq u(x, t)\leq C_2 < \infty$, we have 
$\frac{1}{u^{N-1}(x, t)}-1 \in C^0_{-\tau'}(M)$ as well. We also have $R_{g_0}\in C^{k-2+\alpha}_{-\tau-2}(M)$.
Hence, by the triangle inequality and the product estimate from Lemma \ref{A3},
\begin{align}\notag\displaystyle\frac{R_{g_0}(x)}{u^{N-1}(x, t)} = R_{g_0}(x)\cdot \left(\frac{1}{u^{N-1}(x, t)}-1\right)+ R_{g_0}(x)   \in C^0_{-\tau'-2}(M), \notag
\end{align}
and thus for all $0<\tau_1<\tau'$, $\displaystyle\frac{R_{g_0}(x)}{u^{N-1}(x, t)} \in L^p_{-\tau_1
-2}(M\times [t_0, t_0+1 ]  )$.
%(In fact it is in $C^0_{-2\tau}$ by the product rule of weighted space.)
%\begin{Rem}Do we want to assume $R_0\in C^{k-2+\alpha}_{-2-\tau}(M)$ for $k\geq 2$?\end{Rem}

%Recalling that $\phi(x,t):= u(x, t)-u_\infty(x) \in C^0_{-\tau'}(M)$, $0<C_1\leq u(x, t)\leq C_2 < \infty$ and 
%\eqref{6.7}, 
Putting everything together, we obtain that for any $t_0> 0$ that
\begin{align}\label{6.22}&\|\frac{R_{g_0}\phi}{u^{N-1}}  \|_{L^p_{-\tau_1
-2}(M\times [t_0, t_0+1 ]  )}\\
  \leq& \| \frac{R_{g_0}(x)}{u^{N-1}(x, t)}  \|_{L^p_{-\tau_1
-2}(M\times [t_0, t_0+1] ) }  \cdot 
 \|  \phi(x, t) \|_{C^0_{0}(M )} 
\leq  \frac{C}{t_0^{\delta_0}}. \notag
\end{align}
Now by the weighted Schauder estimate for parabolic equations from \cite[Theorem 5.3]{CZ} (or see Proposition \ref{A4}) applied to \eqref{phieqn}, we have for each $t_0>0$ that
\begin{align}\label{para21p}
&\|   \phi\|_{W^{2,1, p}_{-\tau_1}(M\times [t_0, t_0+1 ])} \notag
\\
\leq& C\left(  \left\|\frac{R_{g_0}\phi}{u^{N-1}} \right\|_{L^p_{-\tau_1-2}(M\times [t_0, t_0+1 ]  )}   + \|  \phi\|_{L^p_{-\tau_1}(M\times [t_0, t_0+1 ]  )} \right),
\end{align}
so substituting \eqref{6.22} and \eqref{6.21} into \eqref{para21p}, we derive
\begin{align}\label{6.19} \|   \phi(x, t)\|_{W^{2,1, p}_{-\tau_1
}(M\times [t_0, t_0+1 ])} \leq  \frac{C}{t_0^{\delta_0}}. \end{align}
Since \eqref{6.19} holds for all $p\geq 1$, we have by Sobolev embedding results for parabolic weighted spaces that (see statement  \eqref{holderembed} of  Lemma \ref{embedlem}),
\begin{align}\notag \|   \phi(x, t)\|_{C^{\alpha,\frac{\alpha}{2}}_{-\tau_1}(M\times [t_0, t_0+1])  } \leq  \frac{C}{t_0^{\delta_0}}. \end{align}
where $C= C(\alpha, \tau_1)$ and we have chosen $p>n+2$ so that $1-\frac{n+2}{p}=\alpha$.

\underline{\bf{Claim 2:}} If \eqref{6.61} is valid for $l\geq 0$, then it holds for $l+2$ when $l\leq k-2$.

% Plugging \eqref{6.21} to the right hand side of Proposition \ref{A4} (One needs to apply 
%Using the same method as in the Claim 1, but with weight $-\tau'$, 
Again from the evolution equation for $\phi$ and the weighted Schauder estimates of Proposition \ref{A4}, we have
\begin{align}\label{6.49}& \|   \phi(x, t)\|_{C^{l+2+\alpha,\frac{l+2+\alpha}{2}}_{-\tau'}(M\times [t_0, t_0+1])
}\notag  \\
\leq &C\left(  \|\frac{R_{g_0}\phi}{u^{N-1}} \|_{C^{l+\alpha,\frac{l+\alpha}{2}}_{-\tau'-2}(M\times [t_0, t_0+1])}   + \|  \phi\|_{C^0_{-\tau'
}(M\times [t_0, t_0+1 ]  )} \right).\end{align}
From \eqref{6.7} we have
\begin{align}\notag \|  \phi\|_{C^0_{-\tau'
}(M\times [t_0, t_0+1 ]  )} \leq  \frac{C}{t_0^{\delta_0}}, \end{align}
So it remains to show that
\begin{align}\label{noloseweight}
  \|\frac{R_{g_0}\phi}{u^{N-1}} \|_{C^{l+\alpha,\frac{l+\alpha}{2}}_{-\tau'-2}(M\times [t_0, t_0+1])}   \leq \frac{C}{t_0^{\delta_0}}.
 \end{align}
 In fact, by the product estimate of Lemma \ref{A3} along with Lemma \ref{6lemv}, and the fact that $u-1\in C^{l+\alpha,\frac{l+\alpha}{2}}_{-\tau'}(M\times [t_0, t_0+1])\hookrightarrow C^{l+\alpha,\frac{l+\alpha}{2}}_{0}(M\times [t_0, t_0+1])$ for $l\leq k-2$,
 \begin{align}\notag
&   \| R_{g_0}(\frac{1}{u^{N-1}}-1) \|_{C^{l+\alpha,\frac{l+\alpha}{2}}_{-\tau'-2}(M\times [t_0, t_0+1])}\\
     \leq&  C \| R_{g_0}\|_{C^{l+\alpha}_{-\tau'-2}(M )}\cdot 
\| u-1\| _{C^{l+\alpha,\frac{l+\alpha}{2}}_{0}(M\times [t_0, t_0+1])}\leq C, \notag
 \end{align}
 and also have $R_{g_0}\in C^{k-2+\alpha}_{-\tau-2}(M)\hookrightarrow C^{l+\alpha}_{-\tau'-2}(M )$ for 
all $l\leq k-2$. 

Therefore 
\begin{align}\notag\frac{R_{g_0}}{u^{N-1}} =R_{g_0} (\frac{1}{u^{N-1}}-1)  + R_{g_0} \in C^{l+\alpha,\frac{l+\alpha}{2}}_{-\tau'-2}(M\times [t_0, t_0+1] ).
 \end{align}
Hence \eqref{noloseweight} follows as
\begin{align}\label{6.51}
 & \|\frac{R_{g_0}\phi}{u^{N-1}} \|_{C^{l+\alpha,\frac{l+\alpha}{2}}_{-\tau'-2}(M\times [t_0, t_0+1])}\notag \\
    \leq &
    \|\frac{R_{g_0}}{u^{N-1}} \|_{C^{l+\alpha,\frac{l+\alpha}{2}}_{-\tau'-2}(M\times [t_0, t_0+1])}\cdot \|   \phi(x, t)\|_{C^{l+\alpha,\frac{l+\alpha}{2}}_{0}(M\times [t_0, t_0+1])}\leq 
   \frac{C}{t_0^{\delta_0}}.
 \end{align}
 Above, the last inequality
 %in \eqref{6.51} 
 uses the inductive hypothesis \eqref{6.61}. By substituting \eqref{6.51} into \eqref{6.49} we have proved Claim 2.
%$$ \|   \phi(x, t)\|_{C^{l+2+\alpha,\frac{l+2+\alpha}{2}}_{-\tau}(M\times [t_0, t_0+1])  }\leq \frac{C}{t_0^{\delta_0}}$$

Claim 1 and 2 together complete the proof of \eqref{6.61}. To conclude, in the definition of the $C^{k+\alpha,\frac{k+\alpha}{2}}_{-\tau}(M\times [t_0, t_0+1])$ spaces as stated in Definition \ref{A1} of the Appendix
%this means \begin{align}\notag \sum_{i+2j\leq k}\sup _{(x,t) \in M\times [t_0, t_0+1 ]} r^{-\tau+i+2j}\left|D_x^i D_t^j \phi(x,t) \right|  \leq \frac{A_k}{t_0^{\delta_0}}. \end{align}
we can take $j=0$ and $t=t_0$, and derive for all $t>0$,
\begin{align}\notag
\|\phi(x, t)  \|_{C_{-\tau'}^{k+\alpha}(M)}  \leq 
\|\phi(x, t)  \|_{C^{k+\alpha,\frac{k+\alpha}{2}}_{-\tau'}(M\times [t_0, t_0+1])}  
%=\sum_{i=0}^{k} \sup _{x\in M} r^{-\tau+i}\left|D_{x}^{i}   \phi(x,t) \right|  
\leq \frac{A_k}{t^{\delta_0}},
\end{align}
and thus \begin{align}\notag u(x,t)-1\rightarrow u_\infty(x)-1  \ \mbox{in} \ C_{-\tau'}^{k+\alpha}(M) \ \mbox{as} \ t\rightarrow \infty.
\end{align}

\end{proof}

\appendix

\section{$C^{k+\alpha}_{-\tau}$ AF metrics along the Yamabe flow}\label{appendixcka}

Here we prove Theorem \ref{AFpreservecka} --- namely, we check that a $C^{k+\alpha}_{-\tau}$ AF metric, $k\geq 2$, continues to be a $C^{k+\alpha}_{-\tau}$ AF metric along a fine solution of the Yamabe flow. First we need to recall the definitions of some parabolic H\"{o}lder spaces from \cite{CZ}.

\begin{Def}[see {\cite[Definition 4.1]{CZ}}]\label{A1}
Let $(M^n,g)$ be a complete Riemannian manifold such that there exists a compact $K\subset M^n$ and a diffeomorphism $\Phi:M^n\backslash K\rightarrow \mb{R}^n\backslash B_{R_0}(0)$, for some $R_0>0$. Let $r\geq 1$ be a smooth function on $M^n$ that agrees under the identification $\Phi$ with the Euclidean radial coordinate $|x|$ in a neighborhood of infinity, and let $\hat{g}$ be a smooth metric on $M^n$ which is equal to the Euclidean metric in a neighborhood of infinity under the identification $\Phi$. Let $\mc{M}=M^n\times[0,T]$. Then with all quantities below computed with respect to the metric $\hat{g}$, we have the following function spaces:

The weighted Lebesgue spaces $L^q_\beta(\mc{M})$, for $q\geq 1$ and weight $\beta\in\mb{R}$, consist of those locally integrable functions on $\mc{M}$ for which the following respective norms are finite:
\begin{align}\notag
\|v\|_{L_{\beta}^{q}(\mc{M})}=\left\{\begin{array}{ll}{\left(\int_0^T\int_{M}|v|^{q} r^{-\beta q-n} d x\right)^{\frac{1}{q}},} & {q<\infty}, \\ {\operatorname{ess} \sup _{\mc{M}}\left(r^{-\beta}|v|\right),} & {q=\infty}.\end{array}\right. 
\end{align}

The weighted Sobolev spaces $W^{k,k/2,q}_\beta(\mc{M})$ are then defined in the usual way with the norms
\begin{align}\notag
\|v\|_{W_{\beta}^{k,k/2, q}(\mc{M})}=\sum_{i+2j\leq k}\left\|D^i_xD^j_t v\right\|_{L_{\beta-i-2j}^{q}(\mc{M})}.
\end{align}

The weighted $C^k$ spaces $C^k_\beta(\mc{M})$ consist of the $C^k$ functions for which the following respective norms are finite:
\begin{align}\notag
\|v\|_{C_{\beta}^{k}(\mc{M})}=\sum_{i+2j\leq k}\sup _{\mc{M}} r^{-\beta+i+2j}\left|D_x^i D_t^j v\right|.
\end{align}

The weighted H\"{o}lder spaces $C_\beta^{k+\alpha,(k+\alpha)/2}(\mc{M})$, $\alpha\in(0,1)$, consist of those $v\in C_\beta^k(\mc{M})$ for which the following respective norms are finite:
\begin{align}\notag
\|v\|_{C_{\beta}^{k+\alpha,(k+\alpha)/2}(\mc{M})}=\|v\|_{C_{\beta}^{k}(\mc{M})}+[v]_{C_\beta^{k+\alpha}(\mc{M})}+\langle v\rangle_{C_\beta^{k+\alpha}(\mc{M})}.
%+\sup _{x \neq y \in M} \min (r(x), r(y))^{-\beta+k+\alpha} \frac{\left|D_{x}^{k} v(x)-D_{x}^{k} v(y)\right|}{|x-y|^{\alpha}}.
\end{align}
Here,
\begin{align}\notag
&[v]_{C_{\beta}^{k+\alpha}\left(\mc{M}\right)}\notag
\\
&=\sum_{i+2 j=k} \sup _{(x, t) \neq(y, s) \in \mc{M}} \min (r(x), r(y))^{-\beta+i+2 j+\alpha} \frac{\left|D_{x}^{i} D_{t}^{j} v(x, t)-D_{x}^{i} D_{t}^{j} v(y, s)\right|}{\delta((x, t),(y, s))^{\alpha}},
\end{align}
where $\delta((x, t),(y, s))=d(x,y)+|t-s|^{\frac{1}{2}}$, and for $k\geq 1$,
\begin{align}\notag
&\langle v\rangle_{C_{\beta}^{k+\alpha}\left(\mc{M}\right)}\notag
\\
&=\sum_{i+2 j=k-1} \sup _{(x, t) \neq(y, s) \in \mc{M}} r(x)^{-\beta+i+2 j+\alpha+1} \frac{\left|D_{x}^{i} D_{t}^{j} v(x, t)-D_{x}^{i} D_{t}^{j} v(x, s)\right|}{|t-s|^{\frac{\alpha+1}{2}}}.
\end{align}
\end{Def}

We now record some embedding results for the spaces above from \cite[Theorem 4.5]{CZ}. Note that \eqref{holderembed} below is a sharpened version of an embedding result which appears in the proof of \cite[Theorem 1.3]{CZ}

\begin{Lem}\label{embedlem}
Under the hypotheses of Definition \ref{A1}, the following inequalities hold:
\begin{enumerate}[(i)]
\item For $1\leq p\leq q\leq\infty$ and $\beta_2<\beta_1$, we have
\begin{align}\notag
\|v\|_{L_{\beta_{1}}^{p}\left(\mc{M}\right)} \leq C\|v\|_{L_{\beta_{2}}^{q}}\left(\mc{M}\right)
\end{align}
\item For $\beta=\beta_1+\beta_2$, $1\leq p,q,s\leq\infty$, and $\frac{1}{p}=\frac{1}{q}+\frac{1}{s}$, we have
\begin{align}\notag
\|v\|_{L_{\beta}^{p}\left(\mc{M}\right)} &\leq\|v\|_{L_{\beta_{1}}^{q}\left(\mc{M}\right)} \mid v \|_{L_{\beta_{2}}^{s}\left(\mc{M}\right)}
\intertext{and}
\|v\|_{C_{\beta}^{\alpha, \alpha / 2}}\left(\mc{M}\right) &\leq\|v\|_{C_{\beta_{1}}^{\alpha, \alpha / 2}\left(\mc{M}\right)} \mid v \|_{C_{\beta_{2}}^{\alpha, \alpha / 2}}\left(\mc{M}\right)^{\prime}
\end{align}
\item\label{holderembed} For $p>n+2$ and $\alpha=1-\frac{n+2}{p}$, we have
\begin{align}\notag
\|v\|_{C^{\alpha,\frac{\alpha}{2}}_{\beta}(\mc{M})}\leq C\|v\|_{W^{2,1,p}_\beta(\mc{M})}.
\end{align}
\end{enumerate}
\end{Lem}

We next state some generalizations of product inequalities for the weighted parabolic Sobolev and H\"{o}lder norms found in \cite[Theorem 4.5]{CZ}. They follow from the definitions of these weighted norms. In the the H\"{o}lder case one must check that $C^{k+\alpha,(k+\alpha)/2}_\beta\hookrightarrow C^{j+\alpha,(j+\alpha)/2}_\beta$ for all $0\leq j\leq k$.

\begin{Lem}\label{A3}
Let $\beta=\beta_1+\beta_2$ with $\beta_1,\beta_2\leq 0$.
\begin{enumerate}[(i)]
\item For $1\leq p,q,s\leq\infty$ with $\frac{1}{p}=\frac{1}{q}+\frac{1}{s}$ and a given nonnegative integer $k$, we have
\begin{align}\notag
\|uv\|_{W^{k,k/2,p}_\beta(\mc{M})}\leq C_1\|u\|_{W^{k,k/2,q}_{\beta_1}(\mc{M})}\|v\|_{W^{k,k/2,s}_{\beta_2}(\mc{M})},
\end{align}
where $C_1$ depends on $\mathcal{M}$ and $k$.

\item For $k\geq 0$, we have
\begin{align}\notag
\|uv\|_{C_\beta^{k+\alpha,(k+\alpha)/2}(\mc{M})}\leq C_2\|u\|_{C_{\beta_1}^{k+\alpha,(k+\alpha)/2}(\mc{M})}\|v\|_{C_{\beta_2}^{k+\alpha,(k+\alpha)/2}(\mc{M})},
\end{align}
where $C_2$ depends on $\mathcal{M}$ and $k$.
\end{enumerate}
\end{Lem}

We now check that an unweighted H\"{o}lder estimate holds for the conformal factor $u$ which evolves along the Yamabe flow.

\begin{Lem}\label{uwholder}
Let $g(t)$ be a fine solution of the Yamabe flow starting from the asymptotically flat manifold $(M^n,g_0)$ on the maximal time interval $[0,T_0)$, given by $g(t)=u(t)^{\frac{4}{n-2}}g_0$ with $u(0)\equiv 1$, and let $T<T_0$. Then given $r_0>0$, there exists a sequence $A_k>0$ for $k=1,2,\ldots$, such that
\begin{align}\notag
\|u\|_{C^{k+\alpha,(k+\alpha)/2}\left(B_{g_{0}}\left(p, r_{0}\right) \times[0, T]\right)} \leq A_k,
\end{align}
independently of the point $p\in M$. Hence there exist uniform bounds on $[0,T]$ for the curvature $Rm(x,t)$ and all of its derivatives.
\end{Lem}
\begin{proof}
By \cite[Theorem 2.4]{CZ}, on a given $[0,T]$ the conformal factor $u(x,t)$ satisfies $0<c_1\leq u(x,t)\leq c_2$ for some constants $c_1,c_2$. Therefore we may apply the Krylov-Safonov estimate for parabolic equations to \eqref{Ypde} and then repeatedly apply the Schauder estimates for parabolic equations (see for instance \cite[Theorem 4.9]{Lieberman}) to obtain the conclusion.
\end{proof}

With the help of the product inequalities of Lemma \ref{A3} we can then adapt the arguments in the proof of \cite[Theorem 5.4]{CZ} to establish higher-order weighted Sobolev and H\"{o}lder estimates for the conformal factor $u(t)$.

\begin{Prop}\label{A4}
Let $u(x,t)$ be a fine solution of the Yamabe flow starting from an asymptotically flat manifold $(M^n,g_0)$ on the maximal time interval $[0,T_0)$, given by $g(t)=u(t)^{\frac{4}{n-2}}g_0$ with $u(0)\equiv 1$. Further let $v=1-u$, and suppose for a fixed $T<T_0$ that $0<\delta\leq u(x,t)\leq C'$ on $[0,T]$.
\begin{enumerate}[(i)]
\item \label{A41}
If $(M^n,g_0)$ is $W^{k,p}_{-\tau}$ asymptotically flat for some $k\geq 2$, then there exists $C=C(n,k,p,\tau,\delta,C')$ such that
\begin{align*}
\|v\|_{W^{k,k/2,p}_{-\tau}(\mc{M})}\leq C\left(\|R_{g_0}\|_{W^{k-2,k/2-1,p}_{-2-\tau}(\mc{M})}+\|v\|_{L^p_{-\tau}(\mc{M})}\right).
\end{align*}

\item\label{A42}
If $(M^n,g_0)$ is $C^{k+\alpha}_{-\tau}$ asymptotically flat for some $k\geq 2$, and we have $\|v\|_{C_0^{k-2+\alpha,(k-2+\alpha)/2}(\mc{M})}\leq C''$, then there exists $C=C(n,k,\tau,\delta,C',C'')$ such that
\begin{align}
\|v\|_{C_{-\tau}^{k+\alpha,(k+\alpha)/2}(\mc{M})}\leq C\left(\|R_{g_0}\|_{C_{-2-\tau}^{k-2+\alpha,(k-2+\alpha)/2}(\mc{M})}+\|v\|_{C_{-\tau}^0(\mc{M})}\right)\label{CZ5.4est}
\end{align}
\end{enumerate}

%Let $u(x,t)$ be a fine solution of the Yamabe flow starting from a $C^{k+\alpha}_{-\tau}$ (for some $k\geq 2$ even) asymptotically flat manifold $(M^n,g_0)$ on the maximal time  interval $[0,T_0)$, given by $g(t)=u(t)^{\frac{4}{n-2}}g_0$ with $u(0)\equiv 1$. Let $v=1-u$, and suppose for a fixed $T<T_0$ that $0<\delta\leq u(x,t)\leq C'$ on $[0,T]$, $\|v\|_{C_0^{k+\alpha,(k+\alpha)/2}(\mc{M})}\leq C''$, and $R_{g_0}\in C^{k-2+\alpha}_{-2-\tau}(M)$. Also let $C'''$ be a constant satisfying $\|v(0)\|_{C^{k+\alpha}_{-\tau}}\leq C'''$. Then there exists a constant $C=C(n,k,\tau,\delta,C',C'',C''')$ such that

\end{Prop}
\begin{proof}
The proofs of \eqref{A41} and \eqref{A42} similar and moreover are straightforward adaptations of the proofs of Theorems 5.3 and 5.4 respectively in \cite{CZ}, once we have the multiplicative inequalities of Lemma \ref{A3}. Indeed the statements of \cite[Theorems 5.3, 5.4]{CZ} are exactly the statements of Proposition \ref{A4} in the case $k=2$. As a result we will only give the proof of \eqref{A42} below, making this choice because in this work we are mainly concerned with $C^{k+\alpha}_{-\tau}$ asymptotically flat manifolds.

Following the proof of \cite[Theorem 5.4]{CZ}, we first use a scaling argument on annuli along with the Schauder estimates for parabolic equations as in \cite[Theorem 4.9]{Lieberman} to obtain
\begin{align}
&\|v\|_{C_{-\tau}^{k+\alpha,(k+\alpha)/2}(E_R\times[0,T])}\notag\label{vinfDelta0}
\\
&\quad\leq C(\|(\partial_t-\Delta_0)v\|_{C_{-2-\tau}^{k-2+\alpha,(k-2+\alpha)/2}(E_R\times[0,T])}+\|v\|_{C^0_{-\tau}(E_R\times[0,T])}),
\end{align}
for $R>R_0$, where $R_0>0$ is as in Definition \ref{funcspaces} and $E_R=\Phi^{-1}(\mb{R}^n\backslash B_R(0))\subset M$. Here $\Delta_0$ denotes the Laplacian with respect to the flat metric defined by $\Phi$ on $E_R$. Define the operator 
\begin{align}\notag
P=h(\Delta_{g_0}-a(n)R_{g_0})=h(g_0^{ij}\frac{\partial^2}{\partial x^i\partial x^j}+b^j\frac{\partial}{\partial x^j}-a(n)R_{g_0}),
\end{align}
where $h=\frac{1}{N(1-v)^{N-1}}$ and in the last expression above we have rewritten $\Delta_{g_0}$ in terms of the Euclidean coordinates given by $\Phi$. We then have
\begin{align}\notag
(\partial_t -P) v=a(n)R_{g_0}.
\end{align}
%from Lemma \ref{uwholder} we have that for fixed $k$, 
We will now compare $(\partial_t-\Delta_0)$ and $(\partial_t- P)$ on $E_R\times[0,T]$. 
%There exist constants $C_1,C_2>0$ so that $\|v\|_{C_0^{k+\alpha,(k+\alpha)/2}}\leq C_1$, $\|h\|_{C_0^{k+\alpha,(k+\alpha)/2}}\leq C_2$. 
Decompose $v=v_0+v_\infty$ using a suitable cutoff function so that $\text{supp}(v_\infty)\subset E_R$; then with the help of the product inequality of Lemma \ref{A3} we can estimate
\begin{align}\label{vinfP}
&\|(\Delta_0-P)v_\infty\|_{C_{-2-\tau}^{k-2+\alpha,(k-2+\alpha)/2}(E_R\times[0,T])}\notag
\\
&\leq \|hg_0^{ij}-\delta_{ij}\|_{C_0^{k-2+\alpha,(k-2+\alpha)/2}(E_R\times[0,T])}\|D^2_x v_\infty\|_{C_{-2-\tau}^{k-2+\alpha,(k-2+\alpha)/2}(E_R\times[0,T])}\notag
\\
&\quad+\|h\|_{C_0^{k-2+\alpha,(k-2+\alpha)/2}(E_R\times[0,T])}\|b\|_{C_{-1}^{k-2+\alpha,(k-2+\alpha)/2}(E_R\times[0,T])}\notag
\\
&\qquad\|D_xv_\infty\|_{C_{-1-\tau}^{k-2+\alpha,(k-2+\alpha)/2}(E_R\times[0,T])}\notag
\\
&\quad+\|h\|_{C_0^{k-2+\alpha,(k-2+\alpha)/2}(E_R\times[0,T])}\|R_{g_0}\|_{C^{k-2+\alpha,(k-2+\alpha)/2}_{-2}(E_R\times[0,T])}\notag
\\
&\qquad\|v_\infty\|_{C^{k-2+\alpha,(k-2+\alpha)/2}_{-\tau}(E_R\times[0,T])}\notag
\\
&\leq C( \|hg_0^{ij}-\delta_{ij}\|_{C_0^{k-2+\alpha,(k-2+\alpha)/2}(E_R\times[0,T])}+\|b\|_{C_{-1}^{k-2+\alpha,(k-2+\alpha)/2}(E_R\times[0,T])}\notag
\\
&\quad+\|R_{g_0}\|_{C^{k-2+\alpha,(k-2+\alpha)/2}_{-2}(E_R\times[0,T])})\|v_\infty\|_{C^{k+\alpha,(k+\alpha)/2}_{-\tau}(E_R\times[0,T])}.
\end{align}
Moreover,
\begin{align}
&\|hg_0^{ij}-\delta_{ij}\|_{C_0^{k-2+\alpha,(k-2+\alpha)/2}(E_R\times[0,T])}+\|b\|_{C_{-1}^{k-2+\alpha,(k-2+\alpha)/2}(E_R\times[0,T])}\notag \label{vinfvanish}
\\
&\quad+\|R_{g_0}\|_{C^{k-2+\alpha,(k-2+\alpha)/2}_{-2}(E_R\times[0,T])}\rightarrow 0
\end{align}
as $R\rightarrow\infty$ because of the asymptotic decay of $v$ and $g_0$. Thus, by taking $R$ sufficiently large we obtain by combining \eqref{vinfDelta0}, \eqref{vinfP}, and \eqref{vinfvanish} that
\begin{align}
\|v_\infty\|_{C_{-\tau}^{k+\alpha,(k+\alpha)/2}(\mc{M})}\leq C(\|(\partial_t-P)v_\infty\|_{C_{-2-\tau}^{k-2+\alpha,(k-2+\alpha)/2}(\mc{M})}+\|v_\infty\|_{C^0_{-\tau}(\mc{M})}.\label{vinfpart}
\end{align}
To deal with the first term on the right, if we let $\zeta_R$ be a suitably chosen cutoff function so that $v_0=\zeta_R v$, $v_\infty=(1-\zeta_R)v$, then we have
\begin{align}\label{inftov}
&\|(\partial_t-P)v_\infty\|_{C^{k-2+\alpha,(k-2+\alpha)/2}_{-2-\tau}} \notag  
\\
&\leq \|(\partial_t-P)v\|_{C^{k-2+\alpha,(k-2+\alpha)/2}_{-2-\tau}}+\|(\partial_t-P)(\zeta_R v)\|_{C^{k-2+\alpha,(k-2+\alpha)/2}_{-2-\tau}}\notag
\\
&\leq 2\|(\partial_t-P)v\|_{C^{k-2+\alpha,(k-2+\alpha)/2}_{-2-\tau}}\notag
\\
&\quad+C_R\|v+|\nabla_{g_0} v| \|_{C^{k-2+\alpha,(k-2+\alpha)/2}_0((E_R\backslash E_{2R})\times[0,T])}.
\end{align}
We leave this for now and turn to consider $v_0$. On the bounded space-time domain $\operatorname{supp}(v_0)\subset\mc{M}$ we can directly apply the Schauder estimates for parabolic equations as in \cite[Theorem 4.9]{Lieberman}, to obtain, similarly to \eqref{vinfpart},
\begin{align}
\|v_0\|_{C_{-\tau}^{k+\alpha,(k+\alpha)/2}(\mc{M})}\leq C(\|(\partial_t-P)v_0\|_{C_{-2-\tau}^{k-2+\alpha,(k-2+\alpha)/2}(\mc{M})}+\|v_0\|_{C^0_{-\tau}(\mc{M})},\label{v0part}
\end{align}
and similar to \eqref{inftov} we also have
\begin{align}\notag
&\|(\partial_t-P)v_0\|_{C^{k-2+\alpha,(k-2+\alpha)/2}_{-2-\tau}}\label{0tov}
\\
&\leq 2\|(\partial_t-P)v\|_{C^{k-2+\alpha,(k-2+\alpha)/2}_{-2-\tau}}\notag
\\
&\quad+C_R\|v+|\nabla_{g_0} v| \|_{C^{k-2+\alpha,(k-2+\alpha)/2}_0((E_R\backslash E_{2R})\times[0,T])}.
\end{align}
Finally, putting \eqref{vinfpart}, \eqref{inftov}, \eqref{v0part}, and \eqref{0tov} together, applying a H\"{o}lder norm interpolation inequality to deal with the $|\nabla g_0 v|$ terms, and recalling that $(\partial_t-P)v=a(n)R_{g_0}$, we obtain the desired estimate \eqref{CZ5.4est}.
\end{proof}

%%%%%%
\begin{comment}
Next we state the results of straightforward computations of the evolution of the scalar curvature $R$ and its derivatives.
\begin{Lem}\label{Revo}
For $k\geq 1$ we have the following differential inequalities under the Yamabe flow $\partial_t g=-R g$:
\begin{align}\notag
\partial_t|\nabla^k R|^2=\Delta|\nabla^k R|^2-2|\nabla^{k+1} R|^2+\sum_{\ell=0}^{k-1}\nabla^\ell{Rm}\ast\nabla^{k-\ell} R\ast\nabla^k R.
\end{align}
\end{Lem}

Next we have a result on the asymptotic decay of derivatives of $R$, assuming decay of derivatives of smaller order of $Rm$.

\begin{Lem}
Let $g(t)$ be a fine solution of the Yamabe flow starting from the asymptotically flat manifold $(M^n,g_0)$ on the maximal time interval $[0,T_0)$, given by $g(t)=u(t)^{\frac{4}{n-2}}g_0$ with $u(0)\equiv 1$. Suppose for a fixed $T<T_0$ that
\begin{align}\notag
|\nabla^\ell Rm|\leq D_\ell r^{-(\tau+\ell+2)},
\end{align}
for all $\ell<k$. Then there exists $C_k>0$ such that
\begin{align}\notag
|\nabla^k R|\leq C_k r^{-(\tau+k+2)}.
\end{align}
\end{Lem}
\begin{proof}
Using Lemma \ref{Revo}, the proof follows the same idea as in \cite[Lemma 5.2]{CZ}.
\end{proof}
\end{comment}
%%%%%

We are now in a position to prove Theorem \ref{AFpreservecka}.

\begin{proof}[Proof of Theorem \ref{AFpreservecka}]
The case $k=2$ is proved in \cite[Theorem 1.3]{CZ}. So we need only concern ourselves with $k\geq 3$, and by induction we may assume that $\|v\|_{C^{k-2+\alpha,(k-2+\alpha)/2}_{-\tau}}$ is bounded on $\mc{M}$.

Clearly by our hypotheses $v$ is bounded in $C^{k-2+\alpha,(k-2+\alpha)/2}_{0}(\mc{M})$. And since $(M^n,g_0)$ is $C^{k+\alpha}_{-\tau}$ asymptotically flat, we also have $R_{g_0}\in C^{k-2+\alpha,(k-2+\alpha)/2}_{-2-\tau}(\mc{M})$. Proposition \ref{A4} therefore implies that $v\in C^{k+\alpha,(k+\alpha)/2}_{-\tau}$, which in turn gives $u(x,t)-1\in C^{k+\alpha}_{-\tau}$ for all $t\in[0,T]$ as desired.

%If $u(x,t)$ on $0\leq t<T_0$ is the conformal factor coresponding to a fine solution of the Yamabe flow on a $C^{k+\alpha}_{-\tau}$ AF manifold $(M^n,g_0)$ with $k\geq 2$ even, then for every fixed $T\in(0,T_0)$ we have $v\in C^0_{-\tau}(M\times[0,T])$ by Theorem \ref{AFpreserve}, and moreover $R_{g_0}\in C^{k-2+\alpha,(k-2+\alpha)/2}_{-2-\tau}(M\times[0,T])$ by the asymptotic flatness of $(M^n,g_0)$ and because $R_{g_0}$ is time-independent. Thus \eqref{CZ5.4est} implies $v\in C^{k+2+\alpha,(k+2+\alpha)/2}_{-\tau}(M\times[0,T])$ so that $v\in C^{k+2+\alpha}_{-\tau}(M)$ along $[0,T]$, which gives the desired conclusion.
\end{proof}

\section{$C^{k+\alpha}_{-\tau}$ conformal deformations}\label{confdefappx}

We will indicate here how the statements Propositions \ref{DMthm} and \ref{CBY} on conformal deformations of AF manifolds with metrics lying in  weighted H\"{o}lder spaces follow from the the corresponding results given in \cite{CB,Maxwell,DM}, which work with weighted Sobolev spaces. We begin by recalling those results.

Corresponding to Proposition \ref{DMthm}, we have:

\begin{Prop}[{\cite[Lemma 4.3]{DM}}]\label{DMappx}
Let $\left(M^{n}, g\right)$ be a $W^{k,p}_{-\tau}$ AF manifold with $k\geq 2$, $k>\frac{n}{p}$, and $\tau\in(0,n-2)$. Suppose $R' \in W_{-2-\tau}^{k-2}$ satisfies $R'\leq R_g$. Then there exists a positive function $\phi$ with $\phi-1\in W^{k,p}_{-\tau}$ such that the scalar curvature of $g'=\phi^{\frac{4}{n-2}}g$ is $R'$. In particular $g'$ is also a $W^{k,p}_{-\tau}$ AF metric.
\end{Prop}

Corresponding to Proposition \ref{CBY} we have:

\begin{Prop}[{\cite[Proposition 3]{DM}}]\label{CBYappx}
Let $(M^n,g)$ be a $W^{k,p}_{-\tau}$ AF manifold with $k\geq 2$, $k>\frac{n}{p}$, and $\tau\in(0,n-2)$. Then the following are equivalent:
\begin{enumerate}[(1)]
    \item We have $Y(M,[g])>0$.
    \item There exists a positive function $\phi$ with $\phi-1\in W^{k,p}_{-\tau}$ such that $\tilde{g}=\phi^{\frac{4}{n-2}}g$ is conformally equivalent to $\tilde{g}$ and $R_{\tilde{g}}\equiv 0$.
\end{enumerate}
\end{Prop}

We note that in \cite{DM}, Proposition \ref{DMappx} is actually stated in the $W^{2,p}_{-\tau}$ setting, but the proof uses the $W^{k,p}_{-\tau}$ results established earlier in \cite{DM}, and the more general $W^{k,p}_{-\tau}$ statement given above follows by the same argument.

Now starting from the hypotheses of Propositions \ref{DMappx} and \ref{CBYappx}, we first observe that $C^{k+\alpha}_{-\tau}\hookrightarrow W^{k,p}_{-\tau'}$ for any $p<\infty$ if $\tau'<\tau$. Therefore Propositions \ref{DMthm} and \ref{CBY} immediately give us the existence of metrics $g'$ and $\tilde{g}$ belonging to $W^{k,p}_{-\tau'}$, and we just need to establish that they additionally belong to $C^{k+\alpha}_{-\tau}$. Moreover, by taking $p$ sufficiently large we have $W^{k,p}_{-\tau'}\hookrightarrow C^\alpha_{-\tau'}$, for any $\tau'<\tau$. We can then conclude the results of both Propositions \ref{DMthm} and \ref{CBY} by proving a regularity estimate arising from this information. Indeed, this is essentially an application of the elliptic theory of weighted H\"{o}lder spaces on punctured regions of $\mathbb{R}^n$ as for instance in \cite{Pacard}.

\begin{Lem}
Let $(M^n,g)$ be a $C^{k+\alpha}_{-\tau'}$ AF manifold, $k\geq 2$, with $\tau\in(0,n-2)$. If $R\in C^{k-2+\alpha}_{-2-\tau}$ and $\phi\in C^{k+\alpha}_{-\tau'}$ for all $\tau'<\tau$ satisfy
\begin{align}
    -a(n)\Delta_{g} \phi+R_g \phi=R \phi^N,\label{Beqn}
\end{align}
then $\phi-1\in C^{k+\alpha}_{-\tau}$.
\end{Lem}
\begin{proof}
Let $\psi=\phi-1$. Elliptic weighted Schauder estimates can be derived by scaling on annuli as in the proof of Proposition \ref{A4} (essentially the same procedure as in \cite{CZ,Bartnik86}), which when applied to \eqref{Beqn} then imply
\begin{align*}
    \|\psi\|_{C^{k+\alpha}_{-\tau'}}\leq C\left(\|R_g\|_{C^{k-2+\alpha}_{-2-\tau'}}+\|R\|_{C^{k-2+\alpha}_{-2-\tau'}}+\|\psi\|_{C^{0}_{-\tau'}}\right).
\end{align*}
Above we have used the product inequalities of Proposition \ref{A3} along with $\psi\in C^\alpha_{0}$ in order to deal with the multiplications of $R_g$ and $R$ against powers of $\psi$. The same estimate would then show us that $\psi\in C^{k+\alpha}_{-\tau}$ once we establish that $\psi\in C^0_{-\tau}$, and we shall now do so.

For ease of notation we assume that we have asymptotic coordinates $\Phi:M\backslash K\rightarrow\mathbb{R}^n\backslash B_1(0)$. Then we can rewrite \eqref{Beqn} in the Euclidean coordinates as
\begin{align}
    -\Delta_0\psi+(\delta^{ij}-g^{ij})D^2_{ij}\psi+b^j D_j\psi+R_g(1+\psi)=R(1+\psi)^N\label{Beqn2}
\end{align}
where $\Delta_0$ is the Euclidean Laplacian and $-\Delta_g=-g^{ij} D^2_{ij}+b^j D_j$. Using that $g\in C^{k+\alpha}_{-\tau}$, $D^2\psi\in C^{k-2+\alpha}_{-2-\tau'}\hookrightarrow C^{0+\alpha}_{-2}$, $D\psi\in C^{k-1+\alpha}_{-1-\tau'}\hookrightarrow C^{0+\alpha}_{-1}$, and $R\in C^{k-2+\alpha}_{-2-\tau}$, all together with the product inequality from Proposition \ref{A3}, we see that \eqref{Beqn2} takes the form
\begin{align}
    -\Delta_0\psi=f\in C^{0+\alpha}_{-2-\tau},\quad\text{on }\mathbb{R}^n\backslash B_1(0).\label{B3eqn}
\end{align}
We now apply the Kelvin transform from $\mathbb{R}^n\backslash B_1(0)$ to $\overline{B_1(0)}\backslash\{0\}$. Following the definitions in \cite{Pacard}, we will denote weighted spaces on $\overline{B_1(0)}\backslash\{0\}$ by for instance ${C^*}^{k+\alpha}_{2-n+\tau}$, in order to distinguish them from weighted spaces on $\mathbb{R}^n\backslash B_1(0)$, which we will still denote by $C^{k+\alpha}_{-\tau}$. We will similarly denote the Kelvin transform of $\psi$ by $\psi^*:\overline{B_1(0)}\backslash\{0\}\rightarrow\mathbb{R}$, and similarly for other functions.

It is then possible to check that $\psi\in C^{k+\alpha}_{-\tau}$ implies $\psi^*\in {C^*}^{k+\alpha}_{2-n+\tau}$. In particular, since $\tau\in(0,n-2)$ we also have $2+n-\tau\in(0,n-2)$, and therefore by \cite[Proposition 2.4]{Pacard}, $\Delta_0:{C^*}^{2+\alpha}_{2-n+\tau,\mathcal{D}}\rightarrow {C^*}^{\alpha}_{2-n+\tau-2}$ is invertible for all $\tau\in(0,n-2)$. Here ${C^*}^{2+\alpha}_{2-n+\tau,\mathcal{D}}$ means those functions of ${C^*}^{2+\alpha}_{2-n+\tau}$ which vanish on $\partial B_1(0)$. Thus, since we already know that $\psi\in C^{k+\alpha}_{-\tau'}$ for any $\tau'<\tau$, we see that $\psi^*\in {C^*}^{2+\alpha}_{2-n+\tau'}$ and
\begin{align}
    -\Delta_0\psi^*&=|x^*|^{-4}f^*\in {C^*}^\alpha_{2-n+\tau'-2},\notag
\intertext{while using \eqref{B3eqn} and \cite[Proposition 2.4]{Pacard}, there exists a $w^*\in {C^*}^{2+\alpha}_{2-n+\tau,\mathcal{D}}$ solving}
    -\Delta_0 w^*&=|x^*|^{-4}f^*\in {C^*}^\alpha_{2-n+\tau-2}.\notag
\end{align}
Hence $\psi^*-w\in {C^*}^{2+\alpha}_{2-n+\tau'}$  must be harmonic in the entire ball $B_1(0)$ since the order $2-n+\tau'>2-n$ implies the singularity at zero can be removed. Thus we see that $\psi^*\in {C^*}^{2+\alpha}_{2-n+\tau}\hookrightarrow {C^*}^0_{2-n_+\tau}$ as well, and transforming back to $\mathbb{R}^n\backslash B_1(0)$ we obtain as desired that $\psi\in C^0_{-\tau}$.
\end{proof}

\section{$L^p$ control of $R$ for small $p$}\label{lowreg}

Here we address regularity issuse that might potentially arise in proving Lemma \ref{monoeps1}, and show that they do not pose any problems.

We first look at a lemma about auxiliary function $w$ constructed in order to deal with the set the scalar curvature vanishes $\{x\in M^n, R(x, t)=0\}$ at time $t$.
%(We can deliberate this motivation later. )
\begin{Lem} \label{lem:4.1}
{\color{red} }
There exists a solution $w$ to $\partial_t w=(n-1)\Delta w+Rw$ with the appropriate spatial decay $w(x, 0)\sim  |x|^{-s}$ and $|\nabla w(x, 0)|\sim  |x|^{-s-1}$ at the initial time and $w(x, 0)>0$ on $M^n$, and $w(x, t)>0$ for all $x\in M^n, t>0$.
\end{Lem}
\begin{proof}[Proof of of Lemma \ref{lem:4.1}] We have already got $L^\infty$ control of $R$ from the Moser iteration, so there exists $\lambda>\sup_{\mathbb{R}^n\times[0,\infty)}|R|$.
By results of \cite{ChauTamYu} the fundamental solution $\Psi(x, t)$ to operator $\partial_t -(n-1)\Delta -R $ exists, and satisfying heat kernel type estimate:
Thus define
$$w(x, t): = \int_M \Psi (x-y, t) w(y, 0) dy$$
satisfying
$\partial_t w=(n-1)\Delta w+Rw$ with initial condition.
We now want to show (1) $w(x, t)>0$ if $w(x, 0)>0$; and (2) $w(x, t)$ decays of order $s$ for all $0<t<t_{max}$ if $w(x, 0)$ decays of order $s$.
%To prove (1), setting $v=e^{\lambda t}w$, we find that
%\begin{equation}\label{eqn:4.8}
%\partial_t v=(n-1)\Delta v+(\lambda+R)v,\end{equation}
%and $(\lambda+R)>0$.
%Moreover, as in Theorem \ref{thm:2.1}, Cheng-Zhu \cite{CZ} proved that asymptotically flatness is preserved under the Yamabe flow \begin{align}\notag \partial_t u^N = \frac{1}{a(n)} \Delta_{g(t)} u- R(g(t)) u \end{align} Thus using the same argument, we can prove for equation (4.8)
It is not hard to see $w$ satisfies conditions of the maximum principle Theorem 4.3 of Ecker and Huisken \cite{EcherHuisken}.
In fact, conditions (i), (ii), (iv) in Theorem 4.3 \cite{EcherHuisken} are obvious. To show condition (iii) $\int_{0, T} \int_M exp(-\alpha_2^2 r^t(p, y)^2) |\nabla w|^2(y)d\mu_t(y) dt<\infty $ for some $\alpha_2>0$, we note that since $\nabla w(x, 0)$ decays, and thus $L^\infty$ on $M^n$,
 $\int_M \Psi (x-y, t) \nabla w(y, 0) d\mu_t(y)$ is $L^\infty$ on $M^n$. Thus
$$\int_M exp(-\alpha_2^2 r^t(p, y)^2) |\nabla w|^2(y)d\mu_t(y) \leq \int_M exp(-\alpha_2^2 r^t(p, y)^2) d\mu_t(y)<\infty .$$

To prove (2) $w(x, t)$ decays of order $s$ if $w(x, 0)$ decays of order $s$, we use the idea of Cheng-Zhu \cite{CZ}, and consider the function
 $f(x, t):= h(x) w(x, t)- C$ where $h(x):= |x|^{s}$, and $C$ is a large enough constant so that $f(x, 0)\leq 0$ using the decay assumption of $w(x, 0)$.
Then $f(x, t)$ satisfies an evolution equation in which its coefficients satisfy all the conditions of the maximum principle of Ecker-Huisken \cite{EcherHuisken} that we used before.
Thus $f(x, t)\leq 0$ as $f(x, 0)\leq 0$. This is to say
$$w(x, t)\leq \frac{C}{|x|^{ s}}.$$
Note $|x|$ is equivalent to the distance function with respect metric $g(t)$, so
$$w(x, t)\leq \frac{C}{d_{g(t)}(x, 0)^{ s}}.$$
\end{proof}

This allows us to prove a version of Lemma \ref{monoeps1} without needing to worry about regularity issues at points when $R=0$. We can control $|R|$ from above by a positive function which satisfies the same decay estimates along with
\begin{align}\notag
    \frac{\partial}{\partial t} A\leq (n-1)\Delta A+RA,
\end{align}
so we can carry out our estimates on $A$ instead.

\begin{proof}[Proof of Lemma \ref{monoeps1}]
%Proposition \ref{prop:4.2}]
Now we can run the rest of the argument to control $L^p$ norms of $R$ when $p<\frac{n}{2}$. Let $A^2:=R^2+w^2$; $\partial_t w=(n-1)\Delta w+Rw$.
By the above explanation, $w>0$. Hence $A>0$.

Note that $w^2$ satisfies $\partial_t w^2=(n-1)\Delta(w^2)-2(n-1)|\nabla w|^2+2R w^2$.  We set the initial condition on $u$ so that it has the same asymptotic decay rate as $R$, ie. $w\sim\epsilon r^{-2-\tau}$. Then by Cauchy-Schwarz
$$|\nabla(R^2+w^2)|^2\leq 4(|R|^2+w^2)(|\nabla|R||^2+|\nabla w|^2).$$ Thus $A^2=R^2+w^2$ satisfies the evolution inequaity

\begin{align}\notag
\partial_t A^2\leq (n-1)\Delta A^2-\frac{n-1}{2}\frac{|\nabla A^2|^2}{A^2}+2RA^2
\end{align}
and therefore (since we have the strict inequality $A>0$),
\begin{align}\notag
\partial_t A\leq (n-1)\Delta A+RA.
\end{align}
The evolution inequality $A$ satisfies is exactly the form of inequality satisfied by $R$, except now we know that $A>|R|\geq 0$.  Thus for any $p$ such that $\int A^p\ dV_t$ is integrable (in particular  $p<\frac{n}{2}$ but close to $\frac{n}{2}$)

\begin{align}\notag
\frac{d}{dt}\int A^p\ dV_t\leq\int pA^{p-1}\partial_t A-\frac{n}{2}RA^p\ dV_t.\end{align}
Hence
\begin{align}\notag
\frac{d}{dt}\int A^p\ dV_t\leq\int pA^{p-1}[(n-1)\Delta A+RA]-\frac{n}{2}RA^p\ dV_t,\end{align}
and
\begin{align}
\frac{d}{dt}\int A^p\ dV_t\leq -p(n-1  )(p-1) \int  A^{p-2} |\nabla A|^2 dV_t+\left(p-\frac{n}{2}\right)\int RA^p\ dV_t.\label{4.2ineq1}
\end{align}

Note $n\geq 3$ implies for $p<\frac{n}{2}$ but close to $\frac{n}{2}$, $p-1>0$. So the first term on the right has a negative sign. In the meanwhile recall for the second term that
\begin{align}
\left |\left(p-\frac{n}{2}\right) \int R A^p \right|\leq\left|p-\frac{n}{2}\right|\|R\|_{L^\frac{n}{2}}\left(\int A^{p\frac{n}{n-2}}\ dV_t\right)^{\frac{n-2}{n}}.\label{4.2ineq2}
\end{align}
Since we have proved $\|R\|_{L^{\frac{n}{2}}}$ is monotonic decreasing, it is bounded depending only on $g_0$, so for $p<\frac{n}{2}$ but very close to $\frac{n}{2}$ this second term can be absorbed by the gradient term, using the Sobolev inequality.  Therefore we see as desired that for such $p$,

\begin{align}\notag
\frac{d}{dt}\int A^p\ dV_t\leq 0.\end{align}

%Therefore $\int A^p\ dV_t$ is bounded for such $p$, and thus
%$\|R\|_{L^p}$ is bounded for such $p$. By interpolation with the $L^\infty$ decay to zero of $|R|$, we see that $\|R\|_{L^{\frac{n}{2}}}$ also goes to zero. Moreover the same argument shows that $\|R\|_{p}$ goes to zero for all $p\in(p_0,\frac{n}{2})$ for some $p\in(1,\frac{n}{2})$.

\end{proof}

%Heuristic argument: (The evolution equation contains a favorable sign when $p$ is slightly smaller than $n/2$)

\printbibliography

%\bibliographystyle{alpha}
%\bibliography{references}

\end{document}